\documentclass[10pt]{amsart}
\usepackage{amsmath}
\usepackage{amscd}
\usepackage{amsxtra}
\usepackage{ifpdf}
\usepackage{amssymb}
\usepackage[mathscr]{eucal}	
\DeclareSymbolFont{pssymbols}     {OMS}{ztmcm}{m}{n}
\DeclareSymbolFontAlphabet{\mathpsscr}   {pssymbols}

\usepackage{amsthm}
\theoremstyle{plain}
\newtheorem{thm}{Theorem}
\newtheorem{cor}[thm]{Corollary}
\newtheorem{prop}[thm]{Proposition}
\newtheorem{lem}[thm]{Lemma}
\newtheorem*{thm*}{Theorem}
\newtheorem*{cor*}{Corollary}
\newtheorem*{prop*}{Proposition}
\newtheorem*{lem*}{Lemma}

\theoremstyle{definition}

\newtheorem*{defn*}{Definition}

\theoremstyle{remark}
\newtheorem{rem}[thm]{Remark}

\newtheorem*{rem*}{Remark}
\newtheorem*{rems*}{Remarks}
\newtheorem*{note*}{Note}

\newcommand{\itemref}[1]{\textup{(\ref{#1})}}
\usepackage[dvips]{graphicx}
\usepackage{psfrag}

\ifpdf                          
  \usepackage[arrow,matrix,tips,frame,curve,pdftex,import]{xy}
\else                           
  \usepackage[arrow,matrix,tips,frame,curve,ps,dvips,import]{xy}
\fi
\SelectTips{cm}{}
\UseTips
\CompileMatrices

\newdir^{ (}{{}*!/-5pt/@^{(}}
\newdir_{ (}{{}*!/-5pt/@_{(}}

\renewcommand{\AA}{{\mathbb A}}
\newcommand{\CC}{{\mathbb C}}
\newcommand{\RR}{{\mathbb R}}
\newcommand{\ZZ}{{\mathbb Z}}
\newcommand{\NN}{{\mathbb N}}
\newcommand{\QQ}{{\mathbb Q}}
\newcommand{\EE}{{\mathbb E}}

\newcommand{\FF}{{\mathbb F}}
\newcommand{\HH}{{\mathbb H}}
\newcommand{\VV}{{\mathbb V}}
\newcommand{\WW}{{\mathbb W}}

\DeclareMathOperator{\codim}{codim}

\DeclareMathOperator{\rank}{rank}

\DeclareMathOperator{\Ker}{Ker}
\renewcommand{\ker}{\Ker}

\DeclareMathOperator{\QQrank}{\QQ-rank}
\DeclareMathOperator{\CCrank}{\CC-rank}
\DeclareMathOperator{\RRrank}{\RR-rank}

\DeclareMathOperator{\supp}{supp}
\newcommand{\geo}{\mathbin{\mathbf o}}

\DeclareMathOperator{\mS}{SS}

\DeclareMathOperator{\pr}{pr}
\DeclareMathOperator{\Cl}{cl}
\newcommand{\cl}[1]{\Cl(#1)}
\renewcommand{\l}{\ell}

\newcommand{\rvvv}[1][]{\rVert_{\if!#1!b\else#1,b\fi}}	
\DeclareMathOperator{\Type}{Type}

\newcommand{\lsb}[1]{{}_{#1}}
\newcommand{\lsp}[1]{{}^{#1}\!}

\newcommand{\tildearrow}{\xrightarrow{\sim}}
\newcommand{\longtildearrow}{\mathrel{\overset{\sim}%
{\longrightarrow}}}
\newbox\arrowbox
\setbox\arrowbox=\hbox{\lower .5ex\rlap{$\longrightarrow$}%
 \raise .5ex\hbox{$\longleftarrow$}}

\newcommand{\G}{\Gamma}

\DeclareMathOperator{\GL}{GL}
\DeclareMathOperator{\SO}{SO}

\DeclareMathOperator{\Hom}{Hom}
\newcommand{\Pl}{\mathscr P}	
\newcommand{\Q}{{\mathscr Q}}	
\newcommand{\R}{{\mathscr R}}
\newcommand{\Center}{Z}
\newcommand{\DerivedGroup}{\mathpsscr D}

\newcommand{\Dbar}{\overline{D}}
\newcommand{\Xbar}{\overline{X}}
\newcommand{\Xhat}{\widehat{X}}
\newcommand{\back}{\backslash}
\newcommand{\Dstar}{D^*}
\newcommand{\Xstar}{X^*}
\newcommand{\Abar}{{\bar{A}}} 

\newcommand{\X}{\mathcal X}

\renewcommand{\i}{i}		
\newcommand{\ihat}{{\hat{\imath}}}
\renewcommand{\j}{j}		

\newcommand{\jhat}{{\hat{\jmath}}}
\DeclareMathOperator{\id}{id}
\newcommand{\nil}{p}	

\newcommand{\wgt}{{h}}    
\newcommand{\cutoff}{\zeta}   

\newcommand{\IpC}{{\mathcal I_p\mathcal C}}

\newcommand{\Sheaf}{\mathcal S}
\newcommand{\A}{\varOmega}
\newcommand{\Asp}{\A_{\,\textup{sp}}}
\newcommand{\Adec}{\widetilde{\A}_{\,\textup{sp}}}

\renewcommand{\sp}{_{\,\textup{sp}}}
\newcommand{\Atilde}{\widetilde{A}}
\newcommand{\inv}{_{\text{inv}}}

\DeclareMathOperator{\Gr}{Gr}

\renewcommand{\L}{\mathscr L}   
\newcommand{\Ltwo}{\A_{(2)}}

\newcommand{\M}{\mathcal M}

\newcommand{\sa}{{\mathfrak a}}

\newcommand{\n}{{\mathfrak n}}

\newcommand{\levi}{{\mathfrak l}}
\newcommand{\p}{{\mathfrak p}}

\newcommand{\al}{\alpha}

\renewcommand{\b}{\beta}
\newcommand{\hsr}{\rho}
\renewcommand{\d}{\delta}

\newcommand{\D}{\Delta}

\newcommand{\g}{\gamma}	
\renewcommand{\u}{\mu}	

\newcommand{\U}{\varUpsilon}

\renewcommand{\o}{\omega}

\makeatletter
\def\sphat{^{\mathchoice{}{}%
 {\,\,\smash[b]{\hbox{\lower4\ex@\hbox{$\m@th\widehat{\null}$}}}}%
 {\,\smash[b]{\hbox{\lower3\ex@\hbox{$\m@th\hat{\null}$}}}}}\,}
\makeatother

\newcommand{\Cat}{\mathscr C}	

\newcommand{\Rep}{\operatorname{\mathfrak M\mathfrak o\mathfrak d}}
\newcommand{\lreg}{{\text{lr}}}

\newcommand{\Derived}{\operatorname{\mathbf D}}

\newcommand{\Complex}{\operatorname{\mathbf C}}

\makeatletter
\def\prime{{\null\prime@\null}}
\mathchardef\prime@="0230
\makeatother

\begin{document}


\author{Leslie Saper}
\address{Department of Mathematics\\ Duke University\\ Box 90320\\ Durham,
NC 27708\\U.S.A.}
\email{saper@math.duke.edu}
\urladdr{http://www.math.duke.edu/faculty/saper}
\title[$L^2$-cohomology I]{\boldmath $L^2$-cohomology of locally symmetric
  spaces, I}
\dedicatory{In memory of Armand Borel}
\thanks{Part of this research was supported in part by the  National Science
Foundation through grants DMS-9870162 and DMS-0502821.
The original manuscript was prepared with the \AmS-\LaTeX\ macro
system and the \Xy-pic\ package.}
\begin{abstract}
Let $X$ be a locally symmetric space associated to a reductive algebraic
group $G$ defined over $\QQ$.  $\L$-modules are a combinatorial analogue of
constructible sheaves on the reductive Borel-Serre compactification
$\Xhat$; they were introduced in \cite{refnSaperLModules}.  That paper also
introduced the micro-support of an $\L$-module, a combinatorial invariant
that to a great extent characterizes the cohomology of the associated
sheaf.  The theory has been successfully applied to solve a number of
problems concerning the intersection cohomology and weighted cohomology of
$\Xhat$ \cite{refnSaperLModules}, as well as the ordinary cohomology of $X$
\cite{refnSaperIHP}.  In this paper we extend the theory so that it covers
$L^2$-cohomology.  In particular we construct an $\L$-module $\Ltwo(X,E)$
whose cohomology is the $L^2$-cohomology $H_{(2)}(X;\EE)$ and we calculate
its micro-support.  As an application we obtain a new proof of the
conjectures of Borel and Zucker.
\end{abstract}
\maketitle
\tableofcontents


\section{Introduction}
\label{sectIntro}
The $L^2$-cohomology $H_{(2)}(X;\EE)$ of an arithmetic locally symmetric
space $X$ plays an important role in geometric analysis and number theory.
In early work, such as \cite{refnBorelStable} and \cite{refnGarlandHsiang},
the application of $L^2$-growth conditions was to single out certain
classes in ordinary cohomology, while later the focus shifted to an
intrinsic notion of $L^2$-cohomology, as in for instance
\cite{refnBorelStableLtwoCohomology}, \cite{refnCheeger}, and
\cite{refnZuckerWarped}.  Zucker conjectured \cite{refnZuckerWarped} that
the $L^2$-cohomology of a Hermitian locally symmetric space $X$ is
isomorphic to the middle-perversity intersection cohomology
$I_pH(\Xstar;\EE)$ of the Baily-Borel-Satake compactification $\Xstar$.
More precisely, the conjecture stated that there is a quasi-isomorphism
$\Ltwo(\Xstar;\EE) \cong \IpC(\Xstar;\EE)$ between complexes of sheaves
which induces the above isomorphism on global cohomology.  Since $\Xstar$
is a projective algebraic variety defined over a number field, the
conjecture is very relevant to Langlands's program and in particular the
study of zeta functions. Zucker \cite{refnZuckerWarped},
\cite{refnZuckerLtwoIHTwo} verified the conjecture in a number of examples.
Borel \cite{refnBorelQRankOne} settled the conjecture in the case where
$\Xstar$ had only one singular stratum; the case of two singular strata was
proved by Borel and Casselman \cite{refnBorelCasselmanQRankTwo}.  The
conjecture in general was resolved in the late 1980's by Stern and the
author \cite{refnSaperSternTwo} and independently by Looijenga
\cite{refnLooijenga}.

{}From the point of view of representation theory, it is natural to consider
the situation where $X$ is an equal-rank locally symmetric space, which is a
more general condition than being Hermitian, and where $\Xstar$ is a Satake
compactification for which all real boundary components of the underlying
symmetric space $D$ are equal-rank.  Borel proposed
\cite[\S6.6]{refnBorelVanishingTheorem}, \cite{refnZuckerLtwoIHTwo} that
Zucker's conjecture be extended to this case.  Soon after
\cite{refnSaperSternTwo} appeared, Stern and the author (unpublished)
verified that their arguments could be extended to settle Borel's
conjecture; this relied partially on a case-by-case analysis.

However for the applications to Langlands's program, one wishes to compute
the local contributions to a fixed-point formula for the action of a
correspondence on $I_pH(\Xstar;\EE)$.  This is complicated by the highly
singular nature of $\Xstar$.  Consequently it is desirable to work on a
less singular compactification of $X$ such as Zucker's reductive
Borel-Serre compactification $\Xhat$ \cite{refnZuckerWarped}, which he
showed \cite{refnZuckerSatakeCompactifications} has a quotient map
$\pi\colon \Xhat\to \Xstar$.  Rapoport \cite{refnRapoportLetterBorel},
\cite{refnRapoport} and independently Goresky and MacPherson
\cite{refnGoreskyMacPhersonWeighted} had conjectured that $I_pH(\Xstar;\EE)
\cong I_pH(\Xhat;\EE)$; more precisely there should be a quasi-isomorphism
$R\pi_*\IpC(\Xhat;\EE) \cong \IpC(\Xstar;\EE)$.  We note also important
related work involving weighted cohomology due to Goresky and MacPherson
and their collaborators \cite{refnGoreskyHarderMacPherson},
\cite{refnGoreskyKottwitzMacPhersonDiscreteSeries},
\cite{refnGoreskyMacPhersonTopologicalTraceFormula}.

Rapoport's conjecture was proved in \cite{refnSaperLModules} for the
equal-rank setting by using the theory of \emph{$\L$-modules} and their
\emph{micro-support}.  An $\L$-module $\M$ is a combinatorial model for a
constructible complex of sheaves on $\Xhat$; the micro-support of an
$\L$-module together with its associated \emph{type} are combinatorial
invariants that to a great extent characterize the cohomology of the
associated sheaf $\Sheaf(\M)$.  The theory is quite general and can be
applied to study many other types of cohomology groups associated to $X$,
for example the weighted cohomology of $\Xhat$ \cite{refnSaperLModules} and
the ordinary cohomology of $X$ \cite{refnSaperIHP}.

Despite the utility of $\L$-modules, they have not yet been used to study
$L^2$-cohomology itself.  (Although $L^2$-cohomology was used as a tool in
\cite{refnSaperLModules} to prove the vanishing theorem recalled in
\S\ref{sectMicroSupport} below, it was not itself the focus of study.)  Of
course, $L^2$-cohomology is by now fairly well-understood; besides the
above references, we note for example other work of Borel and Casselman
\cite{refnBorelCasselman} and Franke \cite{refnFranke}.  Still it would be
valuable to treat $L^2$-cohomology and intersection cohomology within the
same combinatorial framework.  One difficulty that arises is that the
original definition of an $\L$-module does not allow for the infinite
dimensional local cohomology groups which can arise with $L^2$-cohomology.
More seriously, technical analytic problems arise in trying to represent
$L^2$-cohomology as an $\L$-module.

In this paper we overcome these issues and construct a generalized
$\L$-module $\Ltwo(E)$ whose cohomology is the $L^2$-cohomology
$H_{(2)}(X;\EE)$.  We also calculate the micro-support of $\Ltwo(E)$.
These results apply to any locally symmetric space, without the equal-rank
or Hermitian hypothesis.  In a sequel to this paper, we will modify
$\Ltwo(E)$ to obtain an $\L$-module whose cohomology is the ``reduced''
$L^2$-cohomology isomorphic to the space of $L^2$-harmonic differential
forms and compute its micro-support.

As an application of our micro-support calculation and the techniques of
\cite{refnSaperLModules} we obtain here a new proof of the conjectures of
Borel and Zucker.  Elsewhere we will show that a morphism between
$\L$-modules which induces an isomorphism on micro-support and its type
also induces an
isomorphism on global cohomology.  Consequently when the micro-support of
$\Ltwo(E)$ is finite-dimensional (which occurs precisely under the
condition given by Borel and Casselman \cite{refnBorelCasselman}) we
recover Nair's identification of $L^2$-cohomology and weighted cohomology
\cite{refnNair}.  More generally if $(E|_{\lsp0G})^* \cong
\overline{E|_{\lsp0G}}$ then we will establish a relation between the
reduced $L^2$-cohomology, the weighted cohomology, and the intersection
cohomology of $\Xhat$, even beyond the equal-rank situation.  (The
condition $(E|_{\lsp0G})^* \cong \overline{E|_{\lsp0G}}$ is standard in
this context; without it both the $L^2$-cohomology and the weighted
cohomology vanish.)  Unlike the situation of the Borel and Zucker
conjectures, this will not in general be induced from a local isomorphism
on a Satake compactification $\Xstar$.  We note that the relation between
reduced $L^2$-cohomology and weighted cohomology can likely also be proven
using
results of Borel and Garland \cite{refnBorelGarland}, Franke
\cite{refnFranke}, Langlands \cite{refnLanglandsFE}, and Nair
\cite{refnNair}.

The paper begins in \S\ref{sectNotation} by reviewing the notation that we
will use; in particular $D$ will be the symmetric space associated to a
reductive algebraic group $G$ defined over $\QQ$, and $X$ will be the
quotient $\G\back D$ for an arithmetic subgroup $\G\subset G(\QQ)$.  In
\S\ref{sectLTwoCohomology} we briefly recall the definition of
$L^2$-cohomology and the $L^2$-cohomology sheaf.  We give special attention
to the case of a locally symmetric space $X$ with coefficients $\EE$
determined by a regular $G$-module $E$ (that is, where $G \to \GL(E)$ is a
morphism of varieties).  In \S\ref{sectCompactifications} we outline the
construction of the reductive Borel-Serre compactification $\Xhat$ of $X$;
it is a stratified space whose strata are indexed by $\Pl$, the partially
ordered set of $\G$-conjugacy classes of parabolic $\QQ$-subgroups of $G$.
The stratum $X_P$ associated to $P\in \Pl$ is a locally symmetric space
associated to a certain reductive group, namely the Levi quotient
$L_P=P/N_P$, where $N_P$ is the unipotent radical of $P$.

In \S\ref{sectSpecialDifferentialForms} we recall the notion of special
differential forms on $X$ \cite{refnGoreskyHarderMacPherson}; these are
needed in order to associate a sheaf $\Sheaf(\M)$ to an $\L$-module $\M$.
The important fact for us will be that a special differential form on $X$
has a well-defined restriction to a special differential form on any
boundary stratum $X_P$ of $\Xhat$.  The definition of an $\L$-module is
recalled in \S\ref{sectLModules}.  Briefly an $\L$-module $\M$ consists of
a collection of graded regular $L_P$-modules $E_P$, one for each $P\in\Pl$,
together with connecting morphisms $f_{PQ}\colon H(\n_P^Q;E_Q) \to E_P[1]$
whenever $P \le Q$; here $\n_P^Q$ is the Lie algebra of $N_P/N_Q$.  These
data must satisfy a ``differential'' type condition
\eqref{eqnLModuleCondition}.  We also recall the associated sheaf
$\Sheaf(\M)$ on $\Xhat$ as well as pullback and pushforward functors for
$\L$-modules which are analogues of those for sheaves.  In
\S\ref{sectMicroSupport} we recall the micro-support of an $\L$-module and
state a vanishing theorem proved in \cite{refnSaperLModules}.  This theorem
asserts the vanishing of $H(\Xhat;\Sheaf(\M))$ in degrees outside a range
determined by the micro-support of $\M$ and its type.

The new material of the paper begins in \S\ref{sectLocallyRegular}.  The
component $E_P$ of an $\L$-module is actually a complex under the
differential $f_{PP}$; its cohomology represents the local cohomology
$H(\i_P^!\Sheaf(\M))$ with supports along a stratum $\i_P\colon
X_P\hookrightarrow \Xhat$.  Since these groups are often infinite
dimensional for $L^2$-cohomology, we need to generalize the notion of an
$\L$-module to allow $E_P$ to be a locally regular $L_P$-module, that is,
the tensor product of a regular module and a possibly infinite dimensional
vector space on which $L_P$ acts trivially.  We introduce such $\L$-modules
and their associated sheaves in \S\ref{sectLocallyRegular} and verify that
the vanishing theorem continues to hold in this context.

The definition of the $\L$-module $\Ltwo(E)$ is presented in
\S\ref{sectLTwoCohomologyLModule}.  Here is the idea underlying the
definition.  We may assume by induction that $\j_P^*\Ltwo(E)$ has already
been defined, where $j_P\colon U\setminus X_P \hookrightarrow U$ and $U$ is
a neighborhood of some stratum $X_P$.  In order to extend the definition to
all of $U$, one must define a complex $(E_P,f_{PP})$ of locally regular
$L_P$-modules which represents the local $L^2$-cohomology with supports
along $X_P$, together with a map $\bigoplus_{Q>P} f_{PQ}$ from the link
complex $\i_P^*\j_{P*}\j_P^*\Ltwo(E)$ to $(E_P,f_{PP})$.  Zucker's work
\cite{refnZuckerWarped} provides us with a complex of locally regular
$L_P$-modules whose cohomology is the local $L_2$-cohomology along $X_P$
(without supports), namely $(\Ltwo(\Abar_P^G;\HH(\n_P;E),h_P)_\infty,
d_{A_P^G})$; here $\Abar_P^G$ is the compactified split component
transverse to $X_P$, $h_P$ is a certain weight function, and we are taking
germs of forms at infinity.  It is natural to define $(E_P,f_{PP})$ as the
mapping cone (with a degree shift of $-1$) of an attaching map
$\Ltwo(\Abar_P^G;\HH(\n_P;E),h_P)_\infty \to
\i_P^*\j_{P*}\j_P^*\Ltwo(X,E)$.  However the existence of this attaching
map, from forms on $A_P^G$ to forms on smaller split components $A_Q^G$, is
not apparent.  To resolve the problem, we replace
$\i_P^*\j_{P*}\j_P^*\Ltwo(E)$ by a quasi-isomorphic complex of forms on
$A_P^G$, with no additional growth conditions in the new directions, before
forming the mapping cone.

Having defined the $\L$-module $\Ltwo(E)$, we calculate in
\S\ref{sectLocalLtwo} that the associated sheaf $\Sheaf(\Ltwo(E))$ and the
$L^2$-cohomology sheaf $\Ltwo(\Xhat;\EE)$ have the same local cohomology.
However this is not sufficient to establish that they are quasi-isomorphic
since we don't yet know the local quasi-isomorphisms are induced by a
global map of
sheaves.  To construct such a global map requires a complex of forms on $X$
for which both (i) there is a subcomplex whose cohomology is
$L^2$-cohomology, and (ii) there is a restriction map to a similar complex
on any boundary stratum $X_P$.  Special differential forms have the second
property but not the first; smooth forms satisfy the first property but not
the second.  In \S\ref{sectQuasiSpecialDifferentialForms} we introduce the
complex of quasi-special forms and prove it has both desired properties;
this is the technical heart of the paper.  A form is quasi-special if it is
decomposable near any point on the boundary and if the restriction to a
boundary stratum (viewed as a form with coefficients in the sheaf of germs
of forms in the transverse direction) is (recursively) a quasi-special
form.  In \S\ref{sectProofThmLtwoLmodule} we use quasi-special forms to
prove that $\Sheaf(\Ltwo(E))$ and $\Ltwo(\Xhat;\EE)$ are quasi-isomorphic.

Finally the micro-support of $\Ltwo(E)$ is calculated in
\S\ref{sectMicroSupportLtwo} following the analogous calculation for
weighted cohomology in \cite{refnSaperLModules}.  We deduce the conjectures
of Borel and Zucker in \S\ref{sectBorelZuckerConjecture}.

I would like to thank Steve Zucker and Rafe Mazzeo for urging me to write
up this work.  I would also like to thank an anonymous referee for
many thoughtful and insightful comments and suggestions.  I spoke
about these results in July 2004 at the International Conference in Memory
of Armand Borel in Hangzhou.  The $L^2$-cohomology of arithmetic locally
symmetric spaces was a subject that greatly interested Borel, as evidenced
by the many papers he wrote on this subject, particularly during the
1980's.  Thus it seems fitting to dedicate this paper to his memory.

\section{Notation}
\label{sectNotation}
\subsection{Algebraic Groups}
\label{ssectNotationAlgebraicGroups}
For any algebraic group $P$ defined over $\QQ$, let $X(P)$ denote the
regular or rationally defined characters of $P$ and let $X(P)_\QQ$ denote
the subgroup of characters defined over $\QQ$.  Set
\begin{equation*}
\lsp 0 P = \bigcap_{\chi\in X(P)_\QQ} \ker \chi^2.
\end{equation*}
The Lie algebra of $P(\RR)$ will be denoted by $\mathfrak p$.  Let $N_P$
denote the unipotent radical of $P$ and let $L_P = P/N_P$ be its Levi
quotient.  The center of $P$ is denoted by $\Center (P)$ and the derived
group by $\DerivedGroup P$.  Let $S_P$ be the maximal $\QQ$-split torus in
the $\Center (L_P)$ and set $A_P = S_P(\RR)^0$.  We will identify
$X(S_P)\otimes \RR$ with $\sa_P^*$, the dual of the Lie algebra of $A_P$.

Throughout the paper $G$ will be a connected, reductive algebraic group $G$
defined over $\QQ$ and the notation of the previous paragraph will
primarily be applied when $P$ is a parabolic $\QQ$-subgroup of $G$, as we
now assume.  If $P\subseteq Q$ are parabolic $\QQ$-subgroups of $G$, there
are natural inclusions $N_P\subseteq N_Q$ and $A_Q \subseteq A_P$.  We let
$N_P^Q = N_P/N_Q$ denote the unipotent radical of $P/N_Q$ viewed as a
parabolic subgroup of $L_Q$.  There is a natural complement $A_P^Q$ to $A_Q
\subseteq A_P$ which will be recalled in \eqref{eqnSplitDecomposition} and
hence a decomposition $A_P = A_Q \times A_P^Q$.  For $a\in A_P$ we write
$a=a_Qa^Q$ according to this decomposition.  The same notation will be used
for elements of $\sa_P = \sa_Q \oplus \sa_P^Q$ and $\sa_P^* = \sa_Q^*
\oplus \sa_P^{Q*}$.

Let $\D_P\subseteq X(S_P)$ denote the simple weights of the adjoint action
of $S_P$ on the Lie algebra $\n_{P\CC}$ of $N_P$.  (Although this action
depends on the choice of a lift $\widetilde S_P\subseteq P$, its weights do
not.)  By abuse of notation we will call these \emph{roots}.  If $P$ is
minimal, $\D_{P}$ are the simple roots for some ordering of the $\QQ$-root
system of $G$ and we have the \emph{coroots} $\{\al\spcheck\}_{\al\in
\D_{P}}$ in $\sa_{P}$.  In general to define the coroot $\al\spcheck\in
\sa_P$ for $\al \in\D_P$, let $P_0\subset P$ be a minimal parabolic
$\QQ$-subgroup and let $\g$ be the unique element of $\D_{P_0}\setminus
\D_{P_0}^P$ such that $\g|_{\sa_P} = \al$.  Following
\cite{refnArthurTraceFormula} we define $\al\spcheck$ as the projection of
$\g\spcheck \in \sa_{P_0} = \sa_P \oplus \sa_{P_0}^P$ to $\sa_P$.

For parabolic $\QQ$-subgroups $P\subseteq Q$, let $\D_P^Q\subseteq \D_P$
denote those roots which restrict trivially to $A_Q$; they form a basis of
$\sa_P^{Q*}$.  The coroots $\{\al\spcheck\}_{\al\in\D_P^Q}$ are a basis of
$\sa_P^Q$ and we let $\{\b_\al^Q\}_{\al\in\D_P^Q}$ denote the corresponding
dual basis of $\sa_P^{Q*}$.  Likewise let
$\{\b_\al^Q{}\spcheck\}_{\al\in\D_P^Q}$ denote the basis of $\sa_P^Q$ dual
to $\D_P^Q$.  Let 
\begin{gather*}
\sa_P^{Q+} = \{\,H \in \sa_P^Q \mid \langle
\al, H \rangle >0 \text{ for all $\al\in\D_P^Q$}\,\}, \\
\lsp+\sa_P^Q = \{\, H \in \sa_P^Q \mid \langle
\b_\al, H \rangle >0 \text{ for all $\al\in\D_P^Q$}\,\}
\end{gather*}
denote the strictly dominant cone and its open dual cone; similarly define
$\sa_P^{Q*+}$ and $\lsp+\sa_P^{Q*}$.  Set $\sa_P^{+} = \sa_G \oplus
\sa_P^{G+}$, etc.  If $P$ is minimal we may omit it from the
notation.

Let $\cl{Y}$ denote the closure of a subspace $Y$ of a topological space.
We will often use the standard facts that $\al|_{\sa_P^Q} \in -
\cl{\sa_P^{Q*+}}$ for $\al \in \D_P\setminus \D_P^Q$ and that $\sa_P^{Q*+}
\subseteq \lsp+\sa_P^{Q*}$.

Let $\hsr_P\in X(L_P)_\QQ \otimes \QQ$ denote one-half the character by
which $L_P$ acts on $\bigwedge^{\dim \n_P}\n_{P}$; we have $\hsr_P \in
\sa_P^{*+}$.  If $P\subseteq Q$, then $\hsr_P|_{\sa_Q}=\hsr_Q$.  Also
define
\begin{equation}
\label{eqnTauDefinition}
\tau_P^Q = {\sum_{\al \in \D_P^Q} \b_{\al}^Q \in
\sa_P^{Q*+}} \qquad\text{and}\qquad
\tau_P^Q{}\spcheck = {\sum_{\al \in \D_P^Q} \b_{\al}^Q{}\spcheck \in
\sa_P^{Q+}}.
\end{equation}

\subsection{Regular Representations}
\label{ssectNotationRegularRepresentations}
By a \emph{regular representation of $G$} (or a \emph{regular $G$-module})
we mean a finite dimensional complex vector space $E$ together with a
morphism $\sigma\colon G \to \GL(E)$ of algebraic varieties.  In other
words, the representation is rationally defined.  Let $\Rep(G)$ denote the
category of regular $G$-modules.

If $E$ is a regular $G$-module, let $E|_{\lsp0G}$ denote the corresponding
regular $\lsp0G$-module; if $E$ is irreducible or more generally
isotypical, let $\xi_E \in X(S_G)$ denote the character by which $S_G$ acts
on $E$.

If $V$ is an irreducible regular $G$-module, let $E_V$ denote the
$V$-isotypical component, that is, $E_V \cong V \times \Hom_G(V,E)$.

\subsection{Homological Algebra}
\label{ssectNotationHomological}
For an additive category $\Cat$ we let $\Gr(\Cat)$ denote the category of
graded objects of $\Cat$ and we let $\Complex(\Cat)$ denote the category of
(cochain) complexes of objects of $\Cat$.  If $C$ is an object of
$\Gr(\Cat)$ and $k\in\ZZ$, the \emph{shifted object} $C[k]$ is defined by
$C[k]^i = C^{k+i}$.  For a complex $(C,d_C)$ in $\Complex(\Cat)$, define
the \emph{shifted complex} $(C[k],d_{C[k]})$ by $d_{C[k]} = (-1)^k d_C$.
The \emph{mapping cone} $M(f)$ of a morphism $f\colon (C,d_C) \to (D,d_D)$
of complexes is the complex $(C[1]\oplus D, -d_C + d_D + f)$.

Consider a functor $F$ from $\Cat$ to $\Complex(\Cat')$, where $\Cat'$ is
another additive category.  For example, $F$ may be the functor $\EE
\mapsto A(X;\EE)$ sending a local system $\EE$ on a manifold $X$ to the
complex of differential forms with coefficients in $\EE$.  In this case we
extend $F$ to a functor $\Gr(\Cat)\to \Complex(\Cat')$ by defining
\begin{equation}
\label{eqnGradedCoefficients}
F(E) = \bigoplus_k F(E^k)[-k].
\end{equation}
Occasionally we further extend $F$ to a functor $\Complex(\Cat)\to
\Complex(\Cat')$ by means of the associated total complex.

\begin{rem*} In most cases we will make a distinction between a graded
  object $C$ and a complex $(C,d_C)$ created using $C$ and a morphism
  $d_C\colon C \to C[1]$, particularly when working with $\L$-modules.
  This is because often a particular graded object or morphism will enter
  into the definition of several complexes.  However for the complex of
  differential forms we will simply write $A(X;\EE)$ instead of
  $(A(X;\EE),d_X)$ and similarly for the corresponding complex of sheaves.
\end{rem*}

\section{$L^2$-cohomology}
\label{sectLTwoCohomology}
\subsection{\boldmath Definition of $L^2$-cohomology}
\label{ssectLTwoCohomology}
Let $\EE$ be a locally constant sheaf on a manifold $X$, that is, $\EE$ is
the sheaf of locally flat sections of a flat vector bundle on $X$ which we
will also denote $\EE$.  Let $A(X;\EE)$ denote the complex of smooth
differential forms with coefficients in $\EE$; the differential is the
exterior derivative $d=d_X$.  By de Rham's theorem, the cohomology of
$A(X;\EE)$ represents the topological or sheaf cohomology $H(X;\EE)$.
Assume $X$ has a Riemannian metric and $\EE$ has a fiber metric (which may
not be locally constant) and for $\o\in A(X;\EE)$ define the $L^2$-norm
(which may be infinite) by
\begin{equation*}
\smash[t]{\|\o\| = \biggl( \int_X |\o|^2\, dV\biggr)^{\frac12}.}
\end{equation*}
Let $A_{(2)}(X;\EE) \subseteq A(X;\EE)$ denote the subcomplex consisting of
forms $\o$ such that $\o$ and $d\o$ are $L^2$, that is, such that $\|\o\|$,
$\|d\o\|<\infty$.  The cohomology $H_{(2)}(X;\EE)$ of $A_{(2)}(X;\EE)$ is
called the \emph{$L^2$-cohomology of $X$ with coefficients in $\EE$}.  We
also consider the weighted $L^2$-norm%
\footnote{The notation is consistent with \cite{refnFranke} whereas in
  \cite{refnZuckerWarped} our norm would be associated to the weight
  function $\wgt^2$.}
$\|\o\|_\wgt = \|\wgt\o\|$ obtained by multiplying the norm on $\EE$ by a
weight function $\wgt\colon X \to (0,\infty)$.  The cohomology of the
corresponding complex $A_{(2)}(X;\EE,\wgt)$ is the \emph{weighted
$L^2$-cohomology} $H_{(2)}(X;\EE,\wgt)$.  If $X$ is noncompact (our case of
interest) then $H_{(2)}(X;\EE)$ and $H_{(2)}(X;\EE,\wgt)$ are no longer
topological invariants of $X$, but depend on the quasi-isometry class of
$\wgt$ and the metrics.

All of the above extends to the case of a Riemannian orbifold $X$ and a
metrized orbifold locally constant sheaf $\EE$.  The notion of an orbifold
(originally a \emph{$V$-manifold}) was introduced by Satake
\cite{refnSatakeVManifold}; for more details see
\cite{refnChenRuanOrbifoldGromovWitten}.  We also may allow $\EE$ to be
graded (by applying \eqref{eqnGradedCoefficients}).

\subsection{\boldmath Localization of $L^2$-cohomology}
\label{ssectLTwoLocalization}
Let $\A(X;\EE)$ be the complex of sheaves associated to the presheaf
$U\mapsto A(U;\EE)$.  From this point of view, the de Rham isomorphism
follows from the facts that $\A(X;\EE)$ is a fine sheaf and the inclusion
$\EE \to \A(X;\EE)$ is a \emph{quasi-isomorphism} (a morphism which induces
an isomorphism on local cohomology sheaves).  If we apply the analogous
localization to $A_{(2)}(X;\EE)$, the $L^2$ growth conditions disappear and
we obtain the same sheaf $\A(X;\EE)$.  Instead, consider a \emph{partial
compactification} $\Xhat$ of $X$; by this we mean a topological space
$\Xhat$ (not necessarily a manifold) which contains $X$ as a dense
subspace.  Define the \emph{$L^2$-cohomology sheaf} $\Ltwo(\Xhat;\EE)$ to
be the complex of sheaves associated to the presheaf $U\mapsto
A_{(2)}(U\cap X;\EE)$.  If $\Xhat$ is compact and $\Ltwo(\Xhat;\EE)$ is
fine, then the $L^2$-cohomology is isomorphic to the hypercohomology of
$\Ltwo(\Xhat;\EE)$.

\subsection{\boldmath $L^2$-cohomology of Locally Symmetric Spaces}
\label{ssectLocallySymmetricLTwoCohomology}
Let $G$ be a connected reductive algebraic group defined over $\QQ$; we
will use the notation established in \S\ref{ssectNotationAlgebraicGroups}.
Given a maximal compact subgroup $K$ of $G(\RR)$ we obtain a symmetric
space $G(\RR)/KA_G$.   If $K$ and $K'$ are two maximal compact subgroups then
$K' = hKh^{-1}$ for some $h\in \DerivedGroup G(\RR)$ which is unique modulo
$K\cap \DerivedGroup G(\RR)$.  We identify $G(\RR)/K'A_G \tildearrow
G(\RR)/KA_G$ by mapping $gK'A_G \mapsto ghKA_G$; the resulting
$G(\RR)$-homogeneous space is the \emph{symmetric space associated to $G$}
and we denote it $D$.  If $\G \subset G(\QQ)$ is an arithmetic subgroup we
let $X = \G\back D$ denote the corresponding \emph{locally symmetric space
associated to $G$ and $\Gamma$}.

Note that the symmetric space $D$ above may have Euclidean factors since
the maximal $\RR$-split torus $\lsb \RR S_G$ in $\Center(G)$ may be
strictly larger than $S_G$.  Set $\lsb \RR A_G = \lsb \RR S_G(\RR)^0$.  The
choice of a basepoint $x_0\in D$ is equivalent to the choice of a maximal
compact subgroup $K$ and a point $a\in \lsb \RR A_G/A_G$ so that $x_0 =
aKA_G$.  For simplicity we will only consider basepoints with $a=e$.  The
choice of a maximal compact subgroup $K$ in turn determines a unique
involutive automorphism $\theta$ of $G$ (the \emph{Cartan involution})
whose fixed point set in $G(\RR)$ is $K$ \cite[\S1.6]{refnBorelSerre}.
Unless otherwise specified we will not assume that a specific basepoint has
been chosen.

A regular representation $E$ of $G$ determines a locally constant sheaf
$\EE = D \times_\G E$.  In general $X$ is an orbifold and $\EE$ is an
orbifold locally constant sheaf, but we will not mention this explicitly
from now on.  Note that there always exists neat (in particular,
torsion-free) subgroups $\G'\subseteq \G$ with finite index; for such
$\G'$, $\G'\back D$ is smooth and $D \times_{\G'} E$ is an honest flat
vector bundle.

Let $x_0\in D$ be a basepoint and let $KA_G$ and $\theta$ be the associated
stabilizer and Cartan involution.  Choose a Hermitian inner product on $E$
such that $\sigma(g)^* = \sigma(\theta g)^{-1}$ for all $g\in G(\RR)$; such
an inner product always exists and is called \emph{admissible for $x_0$}.
If $E$ is irreducible an admissible inner product is uniquely determined up
to a positive scalar multiple.  The admissible inner product on $E$
determines a fiber metric on $\EE$; in the case that $E$ is isotypical this
is given explicitly as
\begin{equation}
|(gKA_G,v)|_{\EE} = |\xi_E(g)| \cdot |g^{-1}v|_E.
\label{eqnFiberMetric}
\end{equation}
(Properly speaking one should write $|\xi_E^k(g)|^{\frac1k}$ instead of
$|\xi_E(g)|$, where $k\in \NN$ is such that $\xi_E^k\in X(S_G)$ extends to
a character on $G$, but we make this abuse of notation.)  If $x_0' = hx_0$
(where $h\in \DerivedGroup G(\RR)$) is another basepoint then $v \mapsto
|h^{-1}v|_E$ is admissible for $x_0'$; it induces the same fiber metric on
$\EE$.

There exists an invariant nondegenerate bilinear form $B$ on the Lie
algebra $\mathfrak g$ of $G(\RR)$ such that the Hermitian inner product
$\langle X,Y\rangle = B(\overline X,\theta Y)$ is positive definite on
$\mathfrak g_\CC$.  This inner product on $\mathfrak g_\CC$ is admissible
for $x_0$ under the adjoint representation.  In addition it induces an
inner product on $T_{x_0}D$ and hence a $G(\RR)$-invariant Riemannian
metric on $D$.  We give $X$ the induced Riemannian metric.

We now apply \S\ref{ssectLTwoCohomology} to define $A_{(2)}(X;\EE)$ and
$H_{(2)}(X;\EE)$ in this context.  These are well-defined since the choices
above yield quasi-isometric metrics.

\section{Compactifications}
\label{sectCompactifications}
We outline the construction of the Borel-Serre compactification following
\cite{refnBorelSerre} however we use the principal homogeneous spaces
$\mathscr A_P^G$ and $\mathscr N_P(\RR)$ introduced in
\cite{refnSaperLModules} in order to write decompositions independent of a
choice of basepoint.

We also recall the reductive Borel-Serre compactification and use it to
represent $L^2$-cohomology as the hypercohomology of a complex of sheaves.

\subsection{Geodesic Action}
\label{ssectGeodesicAction}
Let $x_0\in D$ be a basepoint with corresponding stabilizer $KA_G$ and
Cartan involution $\theta$.  For $Q$ a parabolic $\QQ$-subgroup of $G$,
there is a unique lift of $L_Q(\RR)$ to $\widetilde L_Q(\RR) \subseteq
Q(\RR)$ which is $\theta$-stable; for $z\in L_Q(\RR)$ let $\tilde z\in
\widetilde L_Q(\RR)$ denote the corresponding lift.  Since $G(\RR) = Q(\RR)
K$, any $x\in D$ may be written as $qKA_G$ for some $q = nr \in Q(\RR) =
N_Q(\RR)\widetilde L_Q(\RR)$.  The \emph{geodesic action} of $z\in
L_Q(\RR)$ on $x \in D$ is defined by
\begin{equation}
z\geo x = n \tilde z  r K A_G.
\end{equation}
For $z=a\in A_Q$ this agrees with the definition given in
\cite[\S3.2]{refnBorelSerre}; in general see
\cite[\S1.1]{refnSaperLModules}.  The geodesic action of $L_Q(\RR)$ is
independent of the choice of $x_0$ and commutes with the action of
$N_Q(\RR)$; the geodesic action of $A_Q$ furthermore commutes with the
action of $Q(\RR)$.

Suppose $P\subseteq Q$ are parabolic $\QQ$-subgroups of $G$.  Since $P/N_Q$
is a parabolic subgroup of $L_Q$, the maximal $\QQ$-split torus in
$\Center(P/N_Q)$ is simply $S_Q$.  Then since $P/N_Q$ projects onto $L_P$,
we may identify $S_Q$ with a subtorus of $S_P$ and $A_Q$ with a subgroup of
$A_P$.  The geodesic action of $a\in A_Q$ is the same whether $a$ is viewed
in $A_Q$ or in $A_P$.

\subsection{Geodesic Decompositions}
\label{ssectGeodesicDecompositions}
We may view $A_G$ as a subgroup of $A_Q$; since $A_G$ acts trivially, the
geodesic action of $A_Q$ descends to $A_Q^G = A_Q/A_G$.  The quotient
$\mathscr A_Q^G = \lsp0Q(\RR)\back D$ is a principal $A_Q^G$-homogeneous
space under the geodesic action and the geodesic quotient $e_Q = A_Q^G
\back D$ is a $\lsp0 Q(\RR)$-homogeneous space.  (A choice of a basepoint
in $D$ determines a basepoint in $\mathscr A_Q^G$ and hence a unique
isomorphism of $A_Q^G$-spaces $\mathscr A_Q^G \cong A_Q^G$ sending the
basepoint to the identity.)  The projections yield
\begin{equation}
D \cong \mathscr A_Q^G \times e_Q,
\label{eqnGeodesicDecomposition}
\end{equation}
an isomorphism of $(A_Q^G\times \lsp0 Q(\RR))$-homogeneous spaces
\cite[\S3.8]{refnBorelSerre}.  (This follows from the identity $Q(\RR) =
\widetilde A_Q \times \lsp0Q(\RR)$ for any lift $\widetilde A_Q$ of $A_Q$.)
We denote by
\begin{equation}
\pr_Q\colon D \longrightarrow \mathscr A_Q^G \qquad \text{and} \qquad
\pr^Q\colon D \longrightarrow e_Q
\end{equation}
the corresponding projections; the latter is called \emph{geodesic
retraction}.  We will propagate this notation and terminology to the
induced decompositions of various quotients and compactifications of $D$ to
be considered below, for example \eqref{eqnScriptADecomposition},
\eqref{eqnCompactifiedGeodesicDecomposition}, and
\eqref{eqnQuotientGeodesicDecomposition}.

For $P\subseteq Q$ note that the geodesic action of $A_P$ on $D$ descends
to an action on $e_Q$.  We now define a subgroup $A_P^Q\subseteq A_P$ which
is complementary to $A_Q \subseteq A_P$ and acts freely on $e_Q$.  Note
there is an injection $X(Q)_\QQ = X(L_Q)_\QQ \hookrightarrow X(P/N_Q)_\QQ =
X(L_P)_\QQ \hookrightarrow X(S_P)$, $\chi \mapsto \chi_P$.  Then set
\begin{equation*}
S_P^Q = \bigl(\bigcap_{\chi\in X(Q)_\QQ} \ker \chi_P\bigr)^0 \subseteq S_P
\end{equation*}
and define $A_P^Q = S_P^Q(\RR)^0$.  There is a direct product
decomposition \cite[1.3(15)]{refnZuckerLtwoIHTwo}%
\footnote{Note $A_P^Q$ is not equal in general to the subgroup $A_{P,Q}$
  defined in \cite{refnBorelSerre} and that the decomposition
  \eqref{eqnSplitDecomposition} is different from $A_P = A_Q \times
  A_{P,Q}$ of \cite[4.3(3)]{refnBorelSerre}.}
\begin{equation}
A_P = A_Q \times A_P^Q
\label{eqnSplitDecomposition}
\end{equation}
and \eqref{eqnGeodesicDecomposition} is an isomorphism of $(A_Q^G\times
A_P^Q)$-homogeneous spaces.

The quotient of \eqref{eqnGeodesicDecomposition} by $\lsp0P(\RR)$ yields an
isomorphism
\begin{equation}
\mathscr A_P^G \cong \mathscr A_Q^G \times \mathscr A_P^Q
\label{eqnScriptADecomposition}
\end{equation}
of $(A_Q^G\times A_P^Q)$-homogeneous spaces, where $\mathscr A_P^Q$ is
defined as $\lsp0P(\RR) \back e_Q = A_Q^G \back \mathscr A_P^G$.  The
quotient of $D \cong \mathscr A_P^G \times e_P$ by $A_Q^G$ yields
\begin{equation}
\label{eqnEQGeodesicDecomposition}
e_Q \cong \mathscr A_P^Q \times e_P.
\end{equation}

\subsection{Partial Compactifications}
There is an isomorphism
\begin{equation*}
A_Q^G \cong (\RR^{>0})^{\D_Q}, \qquad a\longmapsto (a^\al)_{\al\in\D_Q},
\end{equation*}
and we partially compactify by allowing these root coordinates to attain
infinity,
\begin{equation*}
\Abar_Q^G \cong (\RR^{>0}\cup\{\infty\})^{\D_Q}.
\end{equation*}
For all $R\ge Q$, let $o_R\in \Abar_Q^G$ denote the point defined by
\begin{equation*}
o_R^\al = \begin{cases}
  \infty & \text{for $\al\in\D_Q\setminus \D_Q^R$,} \\
  1      & \text{for $\al\in\D_Q^R$.}
	  \end{cases}
\end{equation*}
Then there is a stratification
\begin{equation}
\label{eqnAStratification}
\Abar_Q^G = \coprod_{R\ge Q} A_Q^G\cdot o_R = \coprod_{R\ge Q} A_Q^R\cdot
o_R.
\end{equation}
We sometimes identify $A_Q^R$ with the stratum $A_Q^R\cdot o_R$.

Set $D(Q) = D \times_{A_Q} \Abar_Q^G$; the isomorphism
\eqref{eqnGeodesicDecomposition} extends to
\begin{equation}
D(Q) \cong \bar{\mathscr A}_Q^G \times e_Q,
\label{eqnCompactifiedGeodesicDecomposition}
\end{equation}
where $\bar{\mathscr A}_Q^G = \mathscr A_Q^G \times_{A_Q} \Abar_Q^G$.  The
point $o_Q\in \Abar_Q^G$ determines a well-defined point in $\bar{\mathscr
A}_Q^G$ which we also denote $o_Q$ and in general
\eqref{eqnAStratification} induces a stratification of $\bar{\mathscr
A}_Q^G$.

In general the product decomposition $A_P^G =
A_Q^G \times A_P^Q$ does \emph{not} extend to a product decomposition of
$\Abar_P^G$.%
\footnote{However the product decomposition  $A_P^G = A_Q^G \times
  A_{P,Q}^G$ from \cite[4.3(3)]{refnBorelSerre}
  does extend to $\Abar_P^G$.}
However if $A_P^G(Q) = \{\,a\in \Abar_P^G \mid a^\al<\infty
\text{ for all }\al\in\D_P^Q\,\}$ then \cite[Lemma ~3.6]{refnSaperTilings}
\begin{equation}
\label{eqnCompareAPartialCompactification}
\Abar_Q^G \times A_P^Q \cong A_P^G(Q) \subseteq \Abar_P^G.
\end{equation}
It follows that there is an open inclusion
\begin{equation}
\label{eqnCompareCorners}
\begin{split}
D(Q) = D \times_{A_Q} \Abar_Q^G &= D \times_{A_Q\times A_P^Q} (\Abar_Q^G
\times A_P^Q) \\
&\subseteq D \times_{A_P} \Abar_P = D(P).
\end{split}
\end{equation}
Alternatively, \eqref{eqnScriptADecomposition} and
\eqref{eqnCompareAPartialCompactification} yield
\begin{equation}
\label{eqnCompareScriptAPartialCompactification}
\bar{\mathscr A}_Q^G \times \mathscr A_P^Q \subseteq \bar{\mathscr A}_P^G
\end{equation}
and then by \eqref{eqnEQGeodesicDecomposition} and
\eqref{eqnCompactifiedGeodesicDecomposition} we obtain the inclusion
\begin{equation}
\begin{split}
D(Q) \cong  \bar{\mathscr A}_Q^G \times e_Q &\cong \bar{\mathscr
  A}_Q^G\times \mathscr A_P^Q \times e_P \\
& \subseteq \bar{\mathscr A}_P^G \times e_P \cong D(P).
\end{split}
\end{equation}

\subsection{Borel-Serre Compactification}
Set
\begin{equation}
\Dbar = \bigcup_{Q} D(Q)
\end{equation}
where $Q$ ranges over all parabolic $\QQ$-subgroups of $G$ and we identify
$D(Q)$ with an open subset of $D(P)$ when $P \subseteq Q$.  We identify
$e_Q$ with the subset $\{o_Q\}\times e_Q$ of $D(Q)$ (see
\eqref{eqnCompactifiedGeodesicDecomposition}) and hence obtain a
stratification $\Dbar= \coprod_Q e_Q$.

The group of rational points $G(\QQ)$ acts on $\Dbar$.  The arithmetic
quotient $\Xbar = \G\back \Dbar$ is a compact Hausdorff space called the
\emph{Borel-Serre compactification} of $X$.  The normalizer in $\G$ of a
stratum $e_Q$ of $\Dbar$ is $\G_Q = \G \cap Q$ and the corresponding
stratum of $\Xbar$ is $Y_Q =\G_Q \back e_Q$.  The strata are indexed by the
finite set $\Pl$ of $\G$-conjugacy classes of parabolic $\QQ$-subgroups of
$G$.  To avoid overburdening the notation, we will denote the
$\G$-conjugacy class of $Q$ again simply by $Q$.  With this convention,
$Y_P$ is contained in $\cl{Y_Q}$ if and only if there exists $\gamma \in
\G$ such that $\gamma P \gamma^{-1} \subseteq Q$; in this case we will
write $P \le Q$ and this defines a partial order on $\Pl$.

By reduction theory every point of $e_Q$ has a neighborhood in $\Dbar$ on
which the equivalence relation induced by $\G$ is the same as the
equivalence relation induced by $\G_Q$.  However since $\G_Q$ acts on
\eqref{eqnCompactifiedGeodesicDecomposition} only through the second
factor, we obtain
\begin{equation}
\G_Q\back D(Q) \cong \bar{\mathscr A}_Q^G \times Y_Q.
\label{eqnQuotientGeodesicDecomposition}
\end{equation}
Thus every point $y\in Y_Q$ has a neighborhood in $\Xbar$ (in fact a basis
of neighborhoods) for which \eqref{eqnQuotientGeodesicDecomposition}
induces a decomposition.  Specifically, for $b\in \mathscr A_Q^G$ set
\begin{equation}
\bar{\mathscr A}_Q^G(b) = \cl{\exp(\sa_Q^{G+}) \cdot b}.
\end{equation}
If $O_Q$ is a relatively compact neighborhood of $y$ in $Y_Q$ then for $b$
sufficiently close to $o_Q$,
\begin{equation}
\label{eqnXbarNeighborhood}
\bar{\mathscr A}_Q^G(b) \times O_Q \subseteq \G_Q\back D(Q)
\end{equation}
descends to a neighborhood of $y$ in $\Xbar$.  As $O_Q$ shrinks and $b \to
o_Q$, we obtain a basis of neighborhoods of $y$ in $\Xbar$.

\subsection{Reductive Borel-Serre Compactification}
\label{ssectReductiveBorelSerre}
For a parabolic $\QQ$-subgroup $Q$, the quotient $D_Q = N_Q(\RR)\back e_Q$
is the symmetric space associated to $L_Q$.  On the other hand, the
geodesic action of $L_Q(\RR)$ on $D$ descends to an action on $e_Q$ and the
quotient $\mathscr N_Q(\RR) = \lsp0L_Q(\RR)\back e_Q$ is a principal
$N_Q(\RR)$-homogeneous space.  The projections yield a canonical
decomposition
\begin{equation}
\label{eqnNilpotentDecomposition}
e_Q \cong \mathscr N_Q(\RR) \times D_Q,
\end{equation}
an isomorphism of $N_Q(\QQ) \times \lsp0 L_Q(\RR)$-homogeneous spaces
\cite[\S1.4]{refnSaperLModules}.

The quotient of \eqref{eqnNilpotentDecomposition} by $\G_{N_Q}= \G \cap
N_Q$ yields $\G_{N_Q}\back e_Q \cong \mathscr N_Q(\RR)' \times D_Q$, where
$\mathscr N_Q(\RR)' = \G_{N_Q}\back \mathscr N_Q(\RR)$.  This is a trivial
$\mathscr N_Q(\RR)'$-bundle over $D_Q$.  The further quotient by $\G_{L_Q}
= \G_Q / \G_{N_Q}$ yields a flat $\mathscr N_Q(\RR)'$-bundle
\begin{equation}
\nil_Q\colon Y_Q \longrightarrow X_Q,
\end{equation}
where $X_Q = \G_{L_Q} \back D_Q$ is again a locally symmetric space.  This
is called the \emph{nilmanifold fibration}.

Define $\Xhat = \coprod_{Q\in\Pl} X_Q$ and equip it with the quotient
topology from the natural map $\nil = \coprod_Q \nil_Q \colon \Xbar \to
\Xhat$.  This is the \emph{reductive Borel-Serre compactification} which
was introduced by Zucker \cite{refnZuckerWarped}.  Zucker proves
\cite{refnZuckerSatakeCompactifications} that any Satake compactification
$\Xstar$ \cite{refnSatakeQuotientCompact} arises as a quotient
\begin{equation}
\label{refnSatakeQuotient}
\pi\colon  \Xhat \to \Xstar
\end{equation}
of the reductive Borel-Serre compactification.

The reductive Borel-Serre compactification is stratified by $X_Q$ for $Q\in
\Pl$.  The construction of $\Xhat$ is hereditary in the sense that the
closure of $X_Q$ in $\Xhat$ is itself the reductive Borel-Serre
compactification $\Xhat_Q$ of $X_Q$.  We let $\i_Q$ (resp. $\ihat_Q$)
denote the inclusion map of $X_Q$ (resp. $\Xhat_Q$) into $\Xhat$.

For $x\in X_Q$, let $O_Q$ be a relatively compact neighborhood in $Y_Q$ of
$\nil^{-1}(x)$ which is a union of $\mathscr N_Q(\RR)'$-fibers.  If $b\in
\mathscr A_Q^G$ is sufficiently close to $o_Q$ as in
\eqref{eqnXbarNeighborhood}, the image $V = \nil(\bar{\mathscr A}_Q^G(b)
\times O_Q)$ will be called a \emph{special neighborhood} of $x$ in
$\Xhat$.  Note that
\begin{equation}
\label{eqnSpecialNeighborhood}
V\cap X \cong \mathscr A_Q^G(b) \times O_Q.
\end{equation}
Again as $O_Q$ shrinks and $b \to o_Q$, we obtain a basis of neighborhoods
of $y$ in $\Xhat$.

\subsection{\boldmath $L^2$-cohomology Sheaf on $\Xhat$}
Let $E$ be a regular $G$-module equipped with an admissible inner
product.  Then in \S\ref{ssectLocallySymmetricLTwoCohomology} we have
defined the $L^2$-cohomology $H_{(2)}(X;\EE)$.  By the construction in
\S\ref{ssectLTwoLocalization} we obtain a sheaf $\Ltwo(\Xhat;\EE)$ on the
reductive Borel-Serre compactification.  Zucker \cite{refnZuckerWarped}
shows that this sheaf is fine so its hypercohomology represents
$H_{(2)}(X;\EE)$.

\section{Special Differential Forms}
\label{sectSpecialDifferentialForms}
We recall the the sheaf of \emph{special differential forms}
$\Asp(\Xhat;\EE)$ due to Goresky, Harder, and MacPherson
\cite[(13.2)]{refnGoreskyHarderMacPherson}.  These will be used in the next
section to realize an $\L$-module as a complex of sheaves on $\Xhat$.

Let $\EE$ be a locally constant sheaf on $X$ (for example, $\EE$ could be
induced from a regular representation $E$ of $G$).  For any submanifold
$O_Q\subseteq Y_Q = \G_Q\back e_Q$ which is a union of $\mathscr
N_Q(\RR)'$-fibers, define a form $\o\in A(O_Q;\EE)$ to be
\emph{$N_Q(\RR)$-invariant} if its lift to $e_Q$ is $N_Q(\RR)$-invariant;
denote the subcomplex of such forms by $A\inv(O_Q;\EE)$.

Define the sheaf $\Asp(\Xhat;\EE)$ to have sections over $U\subseteq \Xhat$
consisting of elements $\eta \in A(U\cap X;\EE)$ satisfying the following
condition:
\begin{equation}
\label{eqnSpecial}
\parbox[t]{.9\textwidth}{For every boundary point $x\in U\cap X_Q$, there
exists a special neighborhood $V= \nil(\bar{\mathscr A}_Q^G(b) \times O_Q)
\subseteq U$ of $x$ (see \eqref{eqnSpecialNeighborhood}) such that
$\eta|_{V\cap X} = (\pr^Q)^*\o$ with $\o\in A\inv(O_Q;\EE)$.}
\end{equation}
(Without the $N_Q(\RR)$-invariance condition, these are the forms ``locally
lifted from the boundary'' introduced by Borel
\cite[\S8.1]{refnBorelStable}.)  This condition is stable under exterior
differentiation so we obtain a complex of sheaves $\Asp(\Xhat;\EE)$ which
is fine by \cite[(13.4)]{refnGoreskyHarderMacPherson}.  There are natural
quasi-isomorphisms
\begin{equation}
R\i_{G*}\EE \longtildearrow \Asp(\Xhat;\EE) \longtildearrow \i_{G*} \A(X;\EE)
\label{eqnSpecialFormsResolution}
\end{equation}
induced by inclusion \cite[(13.6)]{refnGoreskyHarderMacPherson}.  The local
normal triviality of $\Xhat$ then implies that $\Asp(\Xhat;\EE)$ is
constructible.  (Recall that a complex of sheaves $\Sheaf$ on $\Xhat$ is
called \emph{constructible} if $H(\i_P^*\Sheaf)$ is locally constant for
all $P\in \Pl$.)

The above construction may be applied to each $\Xhat_R$ for $R\in \Pl$.  We
obtain for each $R$ a functor
\begin{equation}
\label{eqnSpecialDifferentialForms}
\Rep(L_R)\longrightarrow \Complex_\X(\Xhat_R) \qquad \text{given by}
\qquad E_R
\longmapsto \Asp(\Xhat_R;\EE_R),
\end{equation}
where $\Complex_\X(\Xhat_R)$ denotes the category of constructible
complexes of sheaves on $\Xhat_R$.

For $P\le Q \in \Pl$, define the principle $N_P^Q(\RR)$-homogeneous space
$\mathscr N_P^Q(\RR) = N_Q(\RR)\back \mathscr N_P(\RR)$.  Its arithmetic
quotient $\mathscr N_P^Q(\RR)' = (\G_{N_P}/\G_{N_Q})\back \mathscr
N_P^Q(\RR)$ is the fiber of  $\overline X_Q \to
\Xhat_Q$ over $X_P \subseteq \Xhat_Q$.  Let $A\inv(\mathscr
N_P^Q(\RR)';\EE_Q)$ denote the $\G_{L_P}$-module of $N_P(\RR)$-invariant
forms on $\mathscr N_P^Q(\RR)'$ and let $\AA\inv(\mathscr
N_P^Q(\RR)';\EE_Q)$ denote the corresponding locally constant sheaf on
$X_P$.  There
is a well-defined restriction map
\begin{equation*}
\Asp(\Xhat_Q;\EE_Q)\to \ihat_{P*}\Asp(\Xhat_P;\AA\inv(\mathscr
N_P^Q(\RR)';\EE_Q)), \qquad \o \mapsto \o|_{\Xhat_P}.
\end{equation*}
However there is a natural quasi-isomorphism
\cite[(12.15)]{refnGoreskyHarderMacPherson},
\cite[Lem.~4.7]{refnSaperLModules}
\begin{equation}
\label{eqnCoefficientHarmonicProjection}
h_{PQ}\colon \AA\inv(\mathscr N_P^Q(\RR)';\EE_Q) \to \HH(\n_P^Q;E_Q)
\end{equation}
obtained by evaluating a form at a basepoint and applying harmonic
projection.  Thus there is a natural morphism
\cite[Cor.~4.8]{refnSaperLModules}
\begin{equation}
\label{eqnRestrictionSpecialForms}
k_{PQ}:\Asp(\Xhat_Q;\EE_Q) \longrightarrow
\ihat_{P*}\Asp(\Xhat_P;\HH(\n_P^Q;E_Q))
\end{equation}
defined by $k_{PQ}(\o) = \ihat_{P*}\Asp(\Xhat_P;h_{PQ})(\o|_{\Xhat_P})$.  The
morphism $\ihat_P^*(k_{PQ})$ induced on the restriction of these sheaves to
$\Xhat_P$ is a quasi-isomorphism.

Kostant's theorem \cite{refnKostant} implies that there exists a natural
isomorphism
\begin{equation}
\kappa_P^Q\colon H(\n_P^Q;H(\n_Q^R;E_R)) \cong H(\n_P^R;E_R)
\label{eqnNilpotentDegeneration}
\end{equation}
when $P\le Q\le R$ (see for example, \cite[\S0.10.20]{refnSaperLModules});
this isomorphism satisfies
\begin{equation}
\label{eqnNilpotentAssociativity}
\kappa_P^Q\circ H(\n_P^Q; \kappa_Q^R) = \kappa_P^R\circ\kappa_P^Q.
\end{equation}
We will often make use of this isomorphism tacitly.

One may check for $P\le Q\le R$ that
\begin{equation}
\label{eqnComposeRestrictionSpecialForms}
\ihat_{Q*}(k_{PQ})\circ k_{QR} = k_{PR}.
\end{equation}

\section{$\L$-modules}
\label{sectLModules}
An $\L$-module is a combinatorial analogue of a constructible complex of
sheaves on $\Xhat$; it has proved useful in studying various cohomology
groups associated to $\Xhat$.  We recall the definitions following
\cite{refnSaperIHP}, \cite{refnSaperLModules}.

\subsection{\boldmath The Category of $\L$-modules}
Recall that $\Pl$ is the partially ordered finite set of $\G$-conjugacy
classes of parabolic $\QQ$-subgroups of $G$.  Let $\Q\subseteq \Pl$ be a
subset satisfying
\begin{equation}
\text{if  $P\le R \le Q$ where $P$, $Q\in \Q$ and $R\in\Pl$, then
  $R\in\Q$.}
\label{eqnLocallyClosed}
\end{equation}
An \emph{$\L$-module $\M$ on $\Q$} is a pair consisting of
\begin{enumerate}
\item a graded regular $L_P$-module $E_P$ for all $P\in\Q$ and
\item an $L_P$-morphism $f_{PQ}\colon H(\n_P^Q;E_Q)\to E_P[1]$ for all
  $P\le Q\in \Q$
\end{enumerate}
which satisfy the condition that for all $P\le R\in \Q$,
\begin{equation}
\sum_{P\le Q\le R}   f_{PQ}\circ H(\n_P^Q;f_{QR}) = 0.
\label{eqnLModuleCondition}
\end{equation}
Note that as indicated in \S\ref{ssectNotationHomological} we are
implicitly extending the functor $E \mapsto H(n_P^Q;E)$ from $\Rep(L_Q)$ to
$\Gr(\Rep(L_Q))$ by setting $H(\n_P^Q;E_Q) = \bigoplus_k
H(\n_P^Q;E_Q^k)[-k]$.

Let $\M=(E_{\cdot}, f_{\cdot\cdot})$ and $\M'
=(E'_{\cdot}, f'_{\cdot\cdot})$ be $\L$-modules over $\Q$.  A
\emph{morphism} $\phi\colon\M\to \M'$ is a collection of $L_P$-module maps
$\phi_{PQ}\colon H(\n_P^Q;E_Q)\to E'_P$ for all $P\le Q\in \Q$ such that
for all $P\le R\in\Q$,
\begin{equation*}
\sum_{P\le Q\le  R} \phi_{PQ}\circ
H(\n_P^Q;f_{QR}) = \sum_{P\le Q\le R} f'_{PQ}\circ H(\n_P^Q;\phi_{QR}).
\end{equation*}
The composition $\phi'\circ \phi$ of morphisms $\phi:\M\to \M'$ and
$\phi':\M'\to \M''$ is defined by
\begin{equation*}
\smash[t]{
(\phi'\circ\phi)_{PR} = \sum_{P\le Q\le R}
\phi'_{PQ} \circ H(\n_P^Q;\phi_{QR})
}
\end{equation*}
for all $P\le R\in\Q$.  We thus obtain a category $\Rep(\L_\Q)$ of
$\L$-modules over $\Q$.

Given an $\L$-module $\M$ over $\Q$ and $P\in \Q$, equation
\eqref{eqnLModuleCondition} (applied with $P=R$) implies that $f_{PP}\circ
f_{PP}=0$ and thus that $(E_P,f_{PP})$ is a complex.  We obtain a functor
$\M\mapsto \i_P^!\M= (E_P,f_{PP})$ from $\Rep(\L_{\Q})$ to
$\Complex^b(\Rep(L_P))$, the category of bounded complexes of regular
$L_P$-modules.  In the case that $\Q = \{P\}$ this is an equivalence of
categories.  Other complexes contructed from the data of an $\L$-module by
means of \eqref{eqnLModuleCondition} will be introduced in
\S\ref{ssectFunctors}.

\subsection{Realization Functor}
\label{ssectRealization}
A subset $W$ of $\Xhat$ will be called \emph{constructible} if it is a
union of strata.  Since $\Pl$ parametrizes the strata of $\Xhat$, there is
a one-to-one correspondence between constructible subsets of $\Xhat$ and
subsets of $\Pl$, namely $W\mapsto \Pl(W) \equiv \{\,P\in\Pl\mid
X_P\subseteq W\,\}$.  Assume that $W$ is locally closed; this is equivalent
to having $\Pl(W)$ satisfy \eqref{eqnLocallyClosed}.  By abuse of notation,
we will speak of the category of $\L$-modules on $W$ (instead of on
$\Pl(W)$), and write $\Rep(\L_W)$.

For a locally closed constructible set $W$ and $P\in \Pl(W)$, we have
inclusions $\i_P\colon X_P \hookrightarrow W$, $\ihat_P\colon \Xhat_P\cap W
\hookrightarrow W$, $\j_P \colon W \setminus X_P \hookrightarrow W$, and
$\jhat_P \colon W \setminus (\Xhat_P \cap W) \hookrightarrow W$.  Let
$\Derived_\X(W)$ denote the derived category of constructible complexes of
sheaves on $W$.  We define a \emph{realization} functor
\begin{equation*}
\Sheaf_W\colon \Rep(\L_W) \longrightarrow \Derived_\X(W)
\end{equation*}
by setting
\begin{equation}
  \left\{
  \begin{aligned}
    \Sheaf_W(\M) &= \bigoplus_{P\in \Pl(W)} \ihat_{P*} \Asp(\Xhat_P\cap
    W;\EE_P),\\
    d_{\Sheaf_W(\M)} &= \sum_{P\in\Pl(W)} d_P + \sum_{P\le
    Q\in\Pl(W)} \Asp(\Xhat_P\cap W;f_{PQ})\circ k_{PQ},
  \end{aligned}
  \right.
\label{eqnRealizationFunctor}
\end{equation}
where $d_P$ is the differential (exterior differentiation) of
$\Asp(\Xhat_P;\EE_P)$.  The cohomology $H(W;\M)$ of an $\L$-module over $W$
is defined as the hypercohomology $H(W;\Sheaf_W(\M))$.  In fact the sheaf
$\Sheaf_W(\M)$ is fine so $H(W;\M)\cong H(\Sheaf_W(\M)(W))$, the
cohomology of global sections of $\Sheaf_W(\M)$.

\subsection{Functors}
\label{ssectFunctors}
Let $W$ be a locally closed constructible subset of $\Xhat$ and let
$\M\in\Rep(\L_W)$.  For $P\in \Pl(W)$ we define:
\begin{gather}
\i_P^!\M = (E_P,f_{PP}),
\label{eqnLocalCohomologyWithSupports} \\
\i_P^*\M = \Bigl( \bigoplus_{P\le R} H(\n_P^R; E_R), \sum_{P \le R \le S}
  H(\n_P^R; f_{RS}) \Bigr),
\label{eqnLocalCohomology}\\
\i_P^* \j_{P*}\j_P^*\M = \Bigl( \bigoplus_{P<R} H(\n_P^R; E_R), \sum_{P < R
  \le S} H(\n_P^R; f_{RS}) \Bigr).
\label{eqnLinkCohomology}
\end{gather}
(For the individual definitions of $\j_{P*}$ and $\j_P^*$ and more general
functors, see \cite[\S3.4]{refnSaperLModules}.)  By virtue of
\eqref{eqnLModuleCondition} these are complexes and their cohomology groups
are called the \emph{local cohomology supported on $P$}, the \emph{local
cohomology at $P$}, and the \emph{link cohomology at $P$} respectively.
The association $\M\mapsto \i_P^*\M$ is a functor
\begin{equation*}
\Rep(\L_W) \longrightarrow \Rep(\L_{X_P}) \cong \Complex^b(\Rep(L_P))
\end{equation*}
compatible with $\i_P^*$ on the realization, that is,
$\Sheaf_{X_P}(\i_P^*\M) \cong \i_P^*\Sheaf_W(\M)$, and similarly for the
others \cite[\S\S3.4, 4.1]{refnSaperLModules}.

We have short exact sequences corresponding to distinguished triangles in
the derived category
\begin{equation}
\label{eqnStandardShortExactSequence}
0 \longrightarrow \i_P^!\M \longrightarrow \i_P^*\M \longrightarrow \i_P^*
\j_{P*}\j_P^*\M \longrightarrow 0 \ .
\end{equation}
More generally, for $P\le Q \in \Pl(W)$ define
\begin{equation}
\label{eqnLocalCohomologySupportedOnQ}
\i_P^* \ihat_Q^!\M = \Bigl( \bigoplus_{P \le R \le Q} H(\n_P^R; E_R),
\sum_{P \le R \le S \le Q} H(\n_P^R; f_{RS}) \Bigr)
\end{equation}
and 
\begin{equation}
\i_P^* \jhat_{Q*}\jhat_Q^* \M = \Bigl( \bigoplus_{P<R \nleq Q} H(\n_P^R; E_R),
  \sum_{P < R \le S \nleq Q} H(\n_P^R; f_{RS}) \Bigr).
\end{equation}
Again there is a short exact sequence
\begin{equation}
\label{eqnShortExactSequenceWithSupport}
0 \longrightarrow \i_P^*\ihat_Q^!\M \longrightarrow \i_P^*\M
\longrightarrow \i_P^* \jhat_{Q*}\jhat_Q^*\M \longrightarrow 0
\end{equation}
which agrees with \eqref{eqnStandardShortExactSequence} when $Q=P$.

Finally for $P\le Q \in \Pl(W)$ define
\begin{equation}
\label{eqnLocalPQCohomology}
\begin{split}
\i_P^* \i_{Q*}\i_Q^* \M &= \Bigl( \bigoplus_{P \le Q \le R} H(\n_P^R; E_R),
\sum_{P \le Q \le R \le S} H(\n_P^R; f_{RS}) \Bigr) \\
&= H(\n_P^Q;\i_Q^*\M).
\end{split}
\end{equation}
There is a natural surjection
\begin{equation}
\label{eqnPQAttaching}
\i_P^* \M \longrightarrow \i_P^* \i_{Q*}\i_Q^* \M \longrightarrow 0\ .
\end{equation}

\section{Micro-support}
\label{sectMicroSupport}
Let $W \subseteq \Xhat_R$ be an open constructible subset, where $R\in\Pl$.
Let $V$ be an irreducible $L_P$-module for some $P\in \Pl(W)$ and recall
from \S\ref{ssectNotationRegularRepresentations} that $\xi_V$ denotes the
character by which $S_P$ acts on $V$.  Define $Q_V^{\prime R} \ge Q_V^R \ge
P$ such that
\begin{equation}
\label{eqnQQtildeDefinition}
\begin{aligned}
\D_P^{Q_V^R} &= \{\, \al\in\D_P^R\mid \langle \xi_V + \hsr_P, \al\spcheck
\rangle < 0\,\}, \\
\D_P^{Q_V^{\prime R}}& = \{\, \al\in\D_P^R\mid \langle \xi_V + \hsr_P,
\al\spcheck \rangle \le 0\,\}
\end{aligned}
\end{equation}
where $\hsr_P$ is as in \S\ref{ssectNotationAlgebraicGroups}.  Let
$[Q_V^R,Q_V^{\prime R}]$ be the interval $\{\, Q\in\Pl(W) \mid Q_V^R\le
Q\le Q_V^{\prime R}\,\}$ in $\Pl(W)$.

The \emph{micro-support} $\mS(\M)$ of $\M\in\Rep(\L_W)$ is the set of all
irreducible $L_P$-modules $V$ for $P\in \Pl(W)$ satisfying
\begin{enumerate}
\item \label{itemMicroSupportCSC} $(V|_{\lsp0L_P})^* \cong
\overline{V|_{\lsp0L_P}}$, and
\item there exists $Q\in[Q_V^R,Q_V^{\prime R}]$ such that
\begin{equation}
\Hom_{L_P}(V,H(\i_P^*\ihat_Q^!\M)) \neq 0.
\label{eqnMicroSupportNonVanishing}
\end{equation}
\end{enumerate}
The group in \eqref{eqnMicroSupportNonVanishing} is the \emph{type} and is
denoted $\Type_{Q,V}(\M)$.  It is sometimes helpful to also consider the
\emph{weak micro-support} $\mS_w(\M)$ for which condition
\itemref{itemMicroSupportCSC} is omitted.

Now assume $W=\Xhat_R$.  The micro-support of $\M$ and the associated types
control the
non-vanishing of $H(\Xhat_R;\M)$.  Specifically, let $c(V;\M)\le d(V;\M)$
denote
the least and greatest degrees in which $\Type_{Q,V}(\M)$ is nonzero for
any $Q$ as above.  Let $\u$ be the highest weight of $V|_{\lsp0L_P}$ for
a fundamental Cartan subalgebra of $^0\levi_P$ and a choice of
$\theta$-stable ordering.  Let $D_P(V)$ be the symmetric space of the
centralizer in $\lsp0L_P$ of $\u$.  The space $D_P(V)$ depends on the
choice of ordering and we choose the ordering to maximize the dimension of
$D_P(V)$.

\begin{thm}[\protect{\cite[Theorem~10.4]{refnSaperLModules}}]
\label{thmVanishing}
Let $\M$ be an $\L$-module on $\Xhat_R$ where $R\in \Pl$ and set
\begin{align}
\label{eqnLowerBound}
c(\M) &= \inf_{V\in\mS(\M)} \tfrac12(\dim D_P - \dim D_P(V)) +
c(V;\M), \text{ and}\\
\label{eqnUpperBound}
d(\M) &= \sup_{V\in\mS(\M)} \tfrac12(\dim D_P + \dim D_P(V)) +
d(V;\M).
\end{align}
Then $H^i(\Xhat_R;\M)=0$ for $i\notin [c(\M),d(\M)]$.  In particular, if
$\mS(\M)=\emptyset$, then $H(\Xhat_R;\M)=0$.
\end{thm}
We can assume $R=G$.  The proof in \cite{refnSaperLModules} uses the
natural diffeomorphism
$s\colon X \to X_0$ constructed in \cite{refnSaperTilings} onto a
relatively compact open domain $X_0\subseteq X$.  (Although $s$ is
obviously not a retraction, the domain $X_0$ is indeed a deformation
retract of $X$.)  The map $s$ extends to a homeomorphism $\hat s\colon
\Xhat \to \Xhat_0$, where $\Xhat_0$ is constructed analogously to $\Xhat$.
The advantage of $X_0$ over $X$ is that its induced Riemannian metric
extends nondegenerately to the compactification $\Xbar_0$ and similarly for
the various strata of $\Xhat_0$.  Thus we can define a Hodge norm on
sections of $\hat s_*\Sheaf_{\Xhat}(\M)$ and to prove the desired vanishing of
cohomology it suffices to establish an estimate
\begin{equation}
\|d\o\|^2 + \|d^*\o\|^2 \ge C\|\o\|^2
\end{equation}
for global sections $\o$ of $\hat s_*\Sheaf_{\Xhat}(\M)$ in the appropriate
degrees.
A combinatorial generalization of the analytic arguments in
\cite{refnSaperSternTwo} localizes the problem to one for each
$P\in\Pl$.  This local problem is treated by an $L^2$-vanishing theorem
\cite[Theorem~14.1]{refnSaperLModules} originally due to Raghunathan
\cite{refnRaghunathan}, \cite{refnRaghunathanCorrection} (and also
deducible from \cite{refnVoganZuckerman}).

\section{Locally Regular $\L$-modules}
\label{sectLocallyRegular}
We need to generalize slightly the notion of an $\L$-module in order to
apply the theory to $\Ltwo(\Xhat;\EE)$.

A (possibly infinite dimensional) representation $E$ of a connected
reductive $\QQ$-group $L$ is called \emph{locally regular}
\cite[\S1.1.3]{refnGoodmanWallach} if any finite-dimensional subspace $F_0$
is contained within a regular subrepresentation $F\subseteq E$.
Equivalently, $E$ is locally regular if it is the direct limit of its
regular subrepresentations, $E = \varinjlim F$.  Let $\Rep_\lreg(L)$ denote
the category of locally regular $L$-modules.  A locally regular $L$-module
$E$ has a canonical isotypical decomposition $\bigoplus_V V\otimes
\Hom_L(V,E)$, where $V$ ranges over the irreducible regular $L$-modules and
$\Hom_L(V,E)$ is a (possibly infinite dimensional) vector space on which
$L$ acts trivially \cite[\S12.1.1]{refnGoodmanWallach}.

The functor $H(\n_P^Q;\cdot)$ extends to a functor from graded locally
regular $L_Q$-modules to graded locally regular $L_P$-modules.  We can thus
define the category $\Rep_\lreg(\L_W)$ of \emph{locally regular
$\L$-modules over $W$} simply by replacing ``regular'' with ``locally
regular'' in \S\ref{sectLModules}.

The realization functor of a locally regular $\L$-module can be defined as
before provided we specify what is meant by the sheaf of special
differential forms $\Asp(\Xhat_Q;\EE_Q)$ when $E_Q$ is a locally regular
$L_Q$-module.  We define a function from a manifold into the locally
regular representation $E_Q$ to be \emph{smooth} if locally its image lies
in a finite dimensional subspace $F_0$ and the map to $F_0$ is smooth.
This allows one to define $\Asp(\Xhat_Q;\EE_Q)$ as the direct limit of the
special differential forms associated to regular subrepresentations $F_Q$,
\begin{equation}
\label{eqnLimitSpecialDifferentialForms}
\Asp(\Xhat_Q;\EE_Q)  =\varinjlim \Asp(\Xhat_Q;\FF_Q).
\end{equation}

The micro-support $\mS(\M)$ of a locally regular $\L$-module is defined as
before and we have the
\begin{prop}
Theorem ~\textup{\ref{thmVanishing}} remains true for locally regular
$\L$-modules.
\end{prop}

\begin{proof}
Let $\M$ be a locally regular $\L$-module on $W$.  Even though the category
$\Rep_\lreg(\L_W)$ is not abelian it makes sense to write $\widetilde \M
\subseteq \M$ if
\begin{align}
&\widetilde E_P\subseteq E_P, \label{eqnSuperRepresentation}\\
&f_{PQ}(H(\n_P^Q;\widetilde E_Q))\subseteq \widetilde E_P[1]\text{, and}
\label{eqnFPQStable} \\
&\widetilde f_{PQ} = f_{PQ}|_{H(\n_P^Q;\widetilde E_Q)},
\label{eqnFPQRestriction}
\end{align}
for all $P\le Q\in\Pl(W)$.  In this case, the inclusion map defines a
morphism of $\L$-modules.

We will show that
\begin{enumerate}
\item\label{itemLimitCohomology} $H(W;\M) = \varinjlim H(W;\widetilde \M)$, 
\item\label{itemDegreeRanges} $c(\M) = \varinjlim c(\widetilde \M)$, and
  $d(\M) = \varinjlim d(\widetilde \M)$,
\end{enumerate}
where the limits are taken over all regular $\L$-modules
$\widetilde\M\subseteq \M$ and in \itemref{itemDegreeRanges} the limits are
taken in the sense of nets.  If we set $W=\Xhat_R$, this suffices to prove
the proposition since
Theorem ~\ref{thmVanishing} applies to all such $\widetilde\M$.

Let $\phi$ be a global section of $\Sheaf_W(\M)$.  In view of
\eqref{eqnLimitSpecialDifferentialForms}, there exist regular
subrepresentations $F_P\subseteq E_P$ for all $P\in\Pl(W)$ such that $\phi$
is a section of $\bigoplus_P \Asp(\Xhat_P\cap W;\FF_P)$.  It is easy to find
regular representations $\widetilde E_P \supseteq F_P$ for all $P\in\Pl(W)$
satisfying \eqref{eqnSuperRepresentation} and \eqref{eqnFPQStable}.  In
fact, for $P=R$ maximal define $\widetilde E_R = F_R + f_{RR}(F_R)$; in
general, assume that $\widetilde E_Q$ has been defined for $Q>P$ and define
$\widetilde E_P = F_P + f_{PP}(F_P) + \sum_{Q> P}
f_{PQ}(H(\n_P^Q;\widetilde E_Q))$.  Thus there exists a regular $\L$-module
$\widetilde \M\subseteq \M$ such that $\phi$ is a global section of
$\Sheaf_W(\widetilde \M)$, which shows that
\begin{equation*}
\Sheaf_W(\M)(W) = \varinjlim \Sheaf_W(\widetilde \M)(W).
\end{equation*}
Assertion ~\itemref{itemLimitCohomology} follows since cohomology commutes
with direct limits.  To prove assertion ~\itemref{itemDegreeRanges} one can
similarly show for all $P\le Q \in \Pl(W)$ that
\begin{equation*}
H(\i_P^* \ihat_Q^!\M) =  \varinjlim H(\i_P^* \ihat_Q^!\widetilde \M). \qedhere
\end{equation*}
\end{proof}

Henceforth we will refer to a locally regular $\L$-module simply as an
$\L$-module.

\section{The $L^2$-cohomology $\L$-module}
\label{sectLTwoCohomologyLModule}
\subsection{}
\label{ssectPartialLTwoSheaf}
For any $P\in \Pl$, recall from \S\ref{ssectGeodesicDecompositions} that
$\mathscr A_P^G$ denotes the principal $A_P^G$-homogeneous space
$\lsp0P(\RR)\back D$ and let $\i_G\colon \mathscr A_P^G \to \bar{\mathscr
A}_P^G$ denote the inclusion into its partial compactification.  For any
graded complex vector space $E$, let $\EE$ denote the corresponding
constant sheaf on $\mathscr A_P^G$ and let $\A(\mathscr A_P^G;\EE) =
E\otimes \A(\mathscr A_P^G)$ be the complex of sheaves of $\EE$-valued
differential forms (with our usual convention \eqref{eqnGradedCoefficients}
concerning graded coefficients); denote the differential $d_{\mathscr
A_P^G}$.  Consider the pushforward sheaf $\i_{G*}\A(\mathscr A_P^G; \EE)$
on $\bar{\mathscr A}_P^G$ and let $\i_{G*}\A(\mathscr A_P^G; \EE)_\infty$
denote its stalk at $o_P$.  This is likewise a complex.  If $E$ has a
locally regular action of a reductive group $L$, we view
$\i_{G*}\A(\mathscr A_P^G; \EE)_\infty$ as a locally regular $L$-module by
allowing $L$ to act solely on the coefficients.

If $P\le Q \le R$ and $E$ is a locally regular $L_R$-module, we have an
$L_Q$-module $H(\n_Q^R;\i_{G*}\A(\mathscr A_P^G; \EE)_\infty)$.  We view
this as a complex under the differential $H(\n_Q^R;d_{\mathscr A_P^G})$.
On the other hand we also have the complex $\i_{G*}\A(\mathscr A_P^G;
\HH(\n_Q^R; E))_\infty$.  Both of these complexes are naturally isomorphic
as $L_Q$-modules to $ H(\n_Q^R; E) \otimes \i_{G*}\A(\mathscr
A_P^G)_\infty$ but the differentials differ by a sign.  We obtain an
isomorphism of complexes
\begin{equation}
\label{eqnLambdaIsomorphism}
\lambda_Q^R\colon H(\n_Q^R;\i_{G*}\A(\mathscr A_P^G; \EE)_\infty) \to
\i_{G*}\A(\mathscr A_P^G; \HH(\n_Q^R; E))_\infty.
\end{equation}
by multiplying the $H^p(\n_Q^R;E^q)\otimes \i_{G*}\A^r(\mathscr
A_P^G)_\infty$ term by $(-1)^{pr}$.

Assume now that $E$ is a regular $A_P^G$-module.  Given the choice of an
inner product on $E$ and a basepoint $x_0\in {\mathscr A}_P^G$, we define a
fiber metric on $\EE$ and a weight function $\wgt_P$ by
\begin{gather}
\label{eqnCoefficientNorm}
|(a\cdot x_0, v)|_{\EE} = |a^{-1}v|_{E} \qquad (a\in A_P^G),\\
\wgt_P(a\cdot x_0) = \hsr_P(a)^{-1}  \qquad (a\in A_P^G).
\end{gather}
Note the difference of \eqref{eqnCoefficientNorm} from equation
\eqref{eqnFiberMetric} of \S\ref{ssectLocallySymmetricLTwoCohomology}.
When $E$ is isotypical and $A_P^G$ acts via a character $\xi_E$ then
$|(a\cdot x_0, v)|_{\EE} = \xi_E(a)^{-1}|v|_{E}$.  The space $\mathscr
A_P^G$ inherits an invariant Riemannian metric from that of $D$ and hence
we obtain a sheaf of weighted $L^2$-complexes $\Ltwo(\bar{\mathscr A}_P^G;
\EE,\wgt_P)$ as in \S\ref{sectLTwoCohomology}.

More generally, for any $Q\ge P$ define a subcomplex
\begin{equation*}
\A_{(2),Q}(\bar{\mathscr A}_P^G; \EE,\wgt_P) \subseteq \i_{G*}\A(\mathscr
A_P^G; \EE)
\end{equation*}
as the sheaf associated to the subpresheaf
\begin{multline}
\label{eqnLtwoRelative}
U \longmapsto \{\, \psi \mid \psi|_{U\cap (\mathscr A_Q^G\times V)},\
d\psi|_{U\cap (\mathscr A_Q^G\times V)} \in A_{(2)}(U\cap (\mathscr
A_Q^G\times V), \EE,\wgt_P) \\
 \text{ for all $V \subseteq \mathscr A_P^Q$ relatively compact}\,\};
\end{multline}
here we use the decomposition $\mathscr A_P^G \cong \mathscr A_Q^G \times
\mathscr A_P^Q$ from \eqref{eqnScriptADecomposition}.  For $Q=G$ this
equals $\i_{G*}\A(\mathscr A_P^G; \EE)$ while for $Q=P$ we recover
$\Ltwo(\bar{\mathscr A}_P^G; \EE,\wgt_P)$.  Again we let
$\A_{(2),Q}(\bar{\mathscr A}_P^G; \EE,\wgt_P)_\infty$ denote the stalk at
$o_P$.

\subsection{}
Let $E$ be a regular $G$-module with an admissible inner product.  Let
$x_0\in \mathscr A_P^G$ be the projection of the basepoint of $D$ for which
the inner product is admissible.  Apply the considerations of the previous
subsection to the regular $L_P$-module $H(\n_P;E)$ with its induced
admissible inner product.

\begin{lem}
\label{lemPullbackMorphism}
For any $P\le P'\le Q'\le Q$ there exists a natural morphism of complexes
of $L_P$-modules
\begin{equation}
g_{PQ,P'Q'}\colon H(\n_P^{P'}; \A_{(2),Q'}(\bar{\mathscr A}_{P'}^G;
\HH(\n_{P'};E),\wgt_{P'})_\infty)
\longrightarrow \A_{(2),Q}(\bar{\mathscr A}_P^G; \HH(\n_P;E),\wgt_P)_\infty
\end{equation}
which is induced by $\pr_{P'}^*$.  For $P'\le P''\le Q'' \le Q'$ these
morphisms satisfy
\begin{equation}
\label{eqnCompositionCanonicalMorphism}
g_{PQ,P'Q'} \circ H(\n_P^{P'}; g_{P'Q',P''Q''}) = g_{PQ,P''Q''} 
\circ \kappa_P^{P'}.
\end{equation}
\end{lem}

\begin{proof}
The desired morphism will be induced from the diagram
\begin{equation*}
\xymatrix @R-3ex{
{H(\n_P^{P'}; \A_{(2),Q'}(\bar{\mathscr A}_{P'}^G;
  \HH(\n_{P'};E),\wgt_{P'})_\infty)}
\ar@{^{ (}->}[r] \ar@{.>}[dddd]_{g_{PQ,P'Q'}} &
{H(\n_P^{P'}; \i_{G*}\A({\mathscr A}_{P'}^G; \HH(\n_{P'};E))_\infty)}
\ar[dd]^{ \A(\mathscr A_{P'}^G;\kappa_P^{P'})\circ \lambda_P^{P'}}\\
& {} \\
& {\i_{G*}\A({\mathscr A}_{P'}^G; \HH(\n_P;E) )_\infty}
\ar[dd]^{\pr_{P'}^*} \\
& {} \\
{\A_{(2),Q}(\bar{\mathscr A}_P^G; \HH(\n_P;E),\wgt_P)_\infty}
\ar@{^{ (}->}[r] &
{\i_{G*}\A({\mathscr A}_P^G; \HH(\n_P;E))_\infty.}
}
\end{equation*}
Note that although a neighborhood of $o_P$ does not in general decompose
under \eqref{eqnCompareScriptAPartialCompactification}, the inverse image
under $\pr_{P'}$ of a neighborhood of $o_{P'} \in \bar{\mathscr A}_{P'}^G$
does indeed contain a neighborhood of $o_P \in \bar{\mathscr A}_P^G$; thus
$\pr_{P'}^*$ above is well-defined.  All the vertical maps on the
right-hand side are cochain morphisms, so $g_{PQ,P'Q'}$ will be as well.
Equation \eqref{eqnCompositionCanonicalMorphism} is clear from
\eqref{eqnNilpotentAssociativity} 
(note that $g_{PQ,P'Q'}$ involves $\kappa_P^{P'}$).
It remains to check that the $L^2$-condition is preserved.  Since
\eqref{eqnScriptADecomposition} is quasi-isometric to a Riemannian product,
the only issue is in the coefficients.  For elements of
$H(\n_P^{P'};\A_{(2),Q'}(\bar{\mathscr A}_{P'}^G; \WW,\wgt_{P'})_\infty)$,
where $W \subseteq H(\n_{P'};E)$ is an irreducible $L_{P'}$-module, the
coefficients contribute the function $\xi_W(a)^{-2}$ to the $L^2$-integral
by \eqref{eqnCoefficientNorm}.  On the other hand, the image under
$g_{PQ,P'Q'}$ has coefficients in the various irreducible $L_P$-modules $Z
\subseteq H(\n_P^{P'};W)$ which contribute $\xi_Z(a)^{-2}$ to the
$L^2$-integral.  The $L^2$-condition is then preserved since
$\xi_Z|_{A_{P'}^G} = \xi_W$ and $\wgt_P|_{\mathscr A_{P'}^G \times V} \sim
\pr_{P'}^* \wgt_{P'}$ for $V \subseteq \mathscr A_P^{P'}$ relatively
compact.
\end{proof}

\subsection{}
For $\R=\{R_0<\dots < R_k\}\subseteq \Pl$, a totally ordered nonempty
subset, set
\begin{equation}
\label{eqnFRDefinition}
F_\R =  \A_{(2),R_k}(\bar{\mathscr A}_{R_0}^G;
\HH(\n_{R_0};E),\wgt_{R_0})_\infty[-k];
\end{equation}
we view $F_\R$ (which actually only depends on $R_0$ and $R_k$) simply as a
graded locally regular $L_{R_0}$-module, not as a complex.  Define
\begin{equation}
\label{eqnURRDefinition}
u_{\R\R} = (-1)^k d_{\mathscr A_{R_0}^G}\colon F_{\R} \longrightarrow F_{\R}[1]
\end{equation}
and for $\R' = \R \setminus \{R_\l\} \neq \emptyset$, where $0\le \l \le
k$, define
\begin{equation}
\label{eqnURRPrimeDefinition}
  u_{\R\R'}= (-1)^\sigma g_{R_0R_k,R_0'R_{k-1}'}\colon
  H(\n_{R_0}^{R_0'};F_{\R'}) \longrightarrow F_\R[1],
\end{equation}
where
\begin{equation*}
  \sigma = \sigma_{\R\R'} = \smash[t]{
    \begin{cases}
      k      & \text{if $\l=k$ and $R_k=G$,} \\
      \l + 1 & \text{otherwise.} 
  \end{cases}}
\end{equation*}
Define the \emph{$L^2$-cohomology $\L$-module with coefficients $E$} by
\begin{equation}
\Ltwo(E) = 
  \left\{
  \begin{aligned}
    E_P &= \bigoplus_{\substack{\R \\ R_0=P}} F_\R, \\
    f_{PQ} &= \sum_{\substack{\R \supseteq \R' \\ R_0=P, R_0'=Q \\
	\#(\R\setminus \R')\le 1}} u_{\R\R'}.
  \end{aligned}
  \right.
\label{eqnLTwoCohomologyLModule}
\end{equation}
In order to verify condition \eqref{eqnLModuleCondition} is satisfied with
this definition of $f_{PQ}$, we must consider $\R \supseteq \R''$ with
$\#(\R\setminus \R'') \le 2$.  The equations
\begin{equation}
\label{eqnUCondition}
\begin{cases}
u_{\R\R}\circ u_{\R\R} =0 & \text{if  $\#(\R\setminus \R'') = 0$,} \\
u_{\R\R''}\circ H(\n_{R_0}^{R_0''};u_{\R''\R''}) +  u_{\R\R}\circ
u_{\R\R''} = 0 &\text{if  $\#(\R\setminus \R'') = 1$,} \\
\sum_{\R \supset \R' \supset \R''}
u_{\R\R'}\circ H(\n_{R_0}^{R_0'};u_{\R'\R''}) = 0
&\text{if  $\#(\R\setminus \R'') = 2$}
\end{cases}
\end{equation}
follow from $d_{\mathscr A_{R_0}^G}^2=0$, the fact that $g_{PQ,P'Q'}$ is a
morphism of complexes, and \eqref{eqnCompositionCanonicalMorphism}
respectively.  (Note that the last sum in \eqref{eqnUCondition} has only
two terms.)  Condition \eqref{eqnLModuleCondition} follows.

We have used similar notation for the $L^2$-cohomology \emph{$\L$-module}
$\Ltwo(E)$ and the $L^2$-cohomology \emph{sheaf} $\Ltwo(\Xhat;\EE)$.  This
is justified by the
\begin{thm}
\label{thmLtwoLmodule}
There is a natural isomorphism in the derived category
\begin{equation*}
\Sheaf_{\Xhat}(\Ltwo(E)) \cong \Ltwo(\Xhat;\EE).
\end{equation*}
\end{thm}

The proof will appear later in \S\ref{sectProofThmLtwoLmodule}.

\section{Local $L^2$-cohomology}
\label{sectLocalLtwo}
We begin by calculating the local cohomology of the $\L$-module $\Ltwo(E)$
and verifying it agrees with the local cohomology of $\Ltwo(\Xhat;\EE)$.
Strictly speaking this is not needed for the proof of Theorem
~\ref{thmLtwoLmodule}, but it will serve as a model for similar arguments
in \S\ref{sectProofThmLtwoLmodule}.  In addition, the results here will be
used in the calculation of micro-support in \S\ref{sectMicroSupportLtwo}.

\begin{prop}
\label{propLocalLtwo}
The projection $\i_P^* \Ltwo(E) \to \Ltwo(\bar{\mathscr A}_P^G;
\HH(\n_P;E),\wgt_P)_\infty$ is a quasi-isomorphism.
\end{prop}

\begin{proof}
By \eqref{eqnLocalCohomology} and \eqref{eqnLTwoCohomologyLModule} write
\begin{align}
\label{eqnLocalLtwoComplex}
\i_P^* \Ltwo(E) &= \Bigl( \bigoplus_{R\ge P} H(\n_P^R;E_R), \sum_{S\ge R \ge
  P} H(\n_P^R;f_{RS}) \Bigr) \\
&= \Bigl( \bigoplus_{\substack{\R \\  R_0 \ge P}}
H(\n_P^{R_0};F_{\R}),  \sum_{\substack{\R \supseteq \R' \\ R_0' \ge R_0 \ge P
  \\ \#(\R\setminus \R') \le 1}} H(\n_P^{R_0}; u_{\R\R'}) \Bigr)
\notag
\end{align}
where $F_\R$ and $u_{\R\R'}$ are given by
\eqref{eqnFRDefinition}--\eqref{eqnURRPrimeDefinition}.  Define a
decreasing filtration by setting
\begin{equation}
\label{eqnFiltration}
F^p\i_P^* \Ltwo(E) = \bigoplus_{\substack{\R \\  R_0 \ge P \\
    \#(\R\setminus \{P\}) \ge p}}
H(\n_P^{R_0};F_{\R})
\end{equation}
with the induced differential.  The associated graded complex for $p=0$ is
\begin{equation}
\Gr_F^0 \i_P^* \Ltwo(E) = (F_{\{P\}}, d_{\mathscr A_P^G}) =
 \Ltwo(\bar{\mathscr A}_P^G; \HH(\n_P;E),\wgt_P)_\infty,
\end{equation}
while for $p>0$ it is a direct sum
\begin{equation}
\Gr_F^p \i_P^* \Ltwo(E) =
\bigoplus_{\substack{\R \\ R_0 > P \\ \#\R = p}} 
M(-u_{\{P\}\cup\R,\R})[-1]
\end{equation}
of the shifted mapping cones (see \S\ref{ssectNotationHomological})
associated to the chain morphisms
\begin{equation}
-u_{\{P\}\cup\R,\R}\colon (H(\n_P^{R_0};F_{\R}),H(\n_P^{R_0};u_{\R\R}))
\longrightarrow
(F_{\{P\}\cup\R}, u_{\{P\}\cup\R,\{P\}\cup\R})[1].
\end{equation}
Since $u_{\{P\}\cup\R,\R} = \pm g_{PR_{p-1},R_0R_{p-1}}$, Lemma
~\ref{lemPullbackMorphismQuasiIsomorphism} below completes the proof.
\end{proof}

\begin{lem}
\label{lemPullbackMorphismQuasiIsomorphism}
For $P\le P'\le Q$, the morphism $g_{PQ,P'Q}$
\begin{equation}
H(\n_P^{P'}; \A_{(2),Q}(\bar{\mathscr A}_{P'}^G;
\HH(\n_{P'};E), \wgt_{P'})_\infty) \longrightarrow \A_{(2),Q}(\bar{\mathscr
  A}_P^G; \HH(\n_P;E), \wgt_P)_\infty
\end{equation}
from Lemma ~\textup{\ref{lemPullbackMorphism}} is a quasi-isomorphism.
\end{lem}

\begin{proof}
{}From the proof of Lemma ~\ref{lemPullbackMorphism} one sees that
$g_{PQ,P'Q} = \pr_{P'}^* \circ \A(\mathscr A_{P'}^G;\kappa_P^{P'})\circ
\lambda_P^{P'}$, where $\pr_{P'}\colon\mathscr A_{P}^G = \mathscr A_{P'}^G
\times \mathscr A_{P}^{P'} \to \mathscr A_{P'}^G$ and the last two factors
are isomorphisms by \eqref{eqnNilpotentDegeneration} and
\eqref{eqnLambdaIsomorphism}.  Thus it suffices to show
\begin{equation*}
\pr_{P'}^*\colon \A_{(2),Q}(\bar{\mathscr A}_{P'}^G; \HH(\n_{P};E),
\wgt_{P'})_\infty \longrightarrow 
\A_{(2),Q}(\bar{\mathscr A}_P^G; \HH(\n_P;E), \wgt_P)_\infty
\end{equation*}
is a quasi-isomorphism.  Without the $L^2$-conditions this result is
standard: one shows that $\psi \mapsto \psi|_{\mathscr A_{P'}^G \times
\{c\}}$ is a homotopy inverse to $\pr_{P'}^*$ by defining a cochain
homotopy operator that integrates in the $\mathscr A_{P}^{P'}$-factor from
a point $c\in \mathscr A_{P}^{P'}$.  Since the $L^2$-conditions in
\eqref{eqnLtwoRelative} are only imposed on subsets with relatively compact
projection to $\mathscr A_P^{P'}$, these conditions are preserved by the
homotopy operator.
\end{proof}

The homotopy operators used above will be discussed in more detail in
\S\S\ref{ssectBasicHomotopyOperator}---\ref{ssectTruncatedNeighborhoods};
in \S\ref{ssectLTwoHomotopy} a related but more subtle homotopy formula
will be established on $\Ltwo(\bar{\mathscr A}_P^G; \HH(\n_P;E),
\wgt_P)_\infty$.

\begin{cor}
For $P\in \Pl$ and $x\in X_P$, 
\begin{equation*}
H(\i_P^* \Ltwo(E)) \cong H(\Ltwo(\Xhat;\EE)_x).
\end{equation*}
\end{cor}
\begin{proof}
Apply Proposition ~\ref{propLocalLtwo} and Zucker's calculation
\cite[(4.24)]{refnZuckerWarped}. 
\end{proof}

\begin{cor}
\label{corNaturalMorphismLtwo}
For $P\le Q$, the natural morphism \textup(see \eqref{eqnLocalPQCohomology}
and \eqref{eqnPQAttaching}\textup)
\begin{equation}
\label{eqnNaturalMorphism}
H(\i_P^* \Ltwo(E)) \longrightarrow
H(\i_P^*\i_{Q*}\i_Q^* \Ltwo(E))=  H(\n_P^Q;H(\i_Q^* \Ltwo(E)))
\end{equation}
corresponds under Proposition ~\textup{\ref{propLocalLtwo}} to the
composition
\begin{equation}
\label{eqnRetrictionLtwo}
\xymatrix @C-80pt{
{H(\Ltwo(\bar{\mathscr A}_P^G; \HH(\n_P;E),\wgt_P)_\infty)}
\ar[rdd]_{(-1)^\sigma (g_{PQ,PP})_*} & {} & 
{H(\n_P^Q; H(\Ltwo(\bar{\mathscr A}_Q^G; \HH(\n_Q;E),\wgt_Q)_\infty))} \\
{} \\
{} & {H(\A_{(2),Q}(\bar{\mathscr A}_P^G; \HH(\n_P;E),\wgt_P)_\infty)}
\ar[ruu]^{\cong}_*!<5pt,0pt>{\labelstyle(g_{PQ,QQ})_*^{-1}}
}
\end{equation}
where $\sigma=1$ if $Q=G$ and $\sigma=0$ otherwise.
\end{cor}
\begin{proof}
We can assume $P<Q$.  Let $\Phi=(\phi_{\R})_{\R}\in \i_P^* \Ltwo(E)$
represent an element of $H(\i_P^* \Ltwo(E))$; the image $\Phi' \in
\i_P^*\i_{Q*}\i_Q^* \Ltwo(E)$ of $\Phi$ under the natural morphism
\eqref{eqnPQAttaching} is obtained by including only those $\R$ with
$R_0\ge Q$.  The projection of Proposition ~\ref{propLocalLtwo} sends
$\Phi$ to $\phi_{\{P\}}$ and sends $\Phi'$ to $\phi_{\{Q\}}$.  However the
$F_{\{P<Q\}}$-component of the equation $d\Phi=0$ (computed using
\eqref{eqnLocalLtwoComplex}) yields
\begin{equation*}
u_{\{P<Q\},\{P\}}(\phi_{\{P\}}) + u_{\{P<Q\},\{Q\}}(\phi_{\{Q\}}) +
    u_{\{P<Q\},\{P<Q\}}(\phi_{\{P<Q\}})=0.
\end{equation*}
The corollary follows by  applying \eqref{eqnURRDefinition} and
\eqref{eqnURRPrimeDefinition} to yield
\begin{equation*}
(-1)^\sigma g_{PQ,PP}(\phi_{\{P\}}) - g_{PQ,QQ}(\phi_{\{Q\}}) -
    d_{\mathscr A_P^G}(\phi_{\{P<Q\}})=0. \qedhere
\end{equation*}
\end{proof}

\section{Quasi-special Differential Forms}
\label{sectQuasiSpecialDifferentialForms}
Proposition ~\ref{propLocalLtwo} of the preceding section does not imply
Theorem ~\ref{thmLtwoLmodule} of \S\ref{sectLTwoCohomologyLModule} since it
only provides a local quasi-isomorphism which may not arise from a global
morphism $\Sheaf_{\Xhat}(\Ltwo(E)) \to \Ltwo(\Xhat;\EE)$.  A global
quasi-isomorphism, albeit in the opposite direction, will be constructed in
\S\ref{sectProofThmLtwoLmodule}.  However this morphism can only be defined
if we replace $\Ltwo(\Xhat;\EE)$ by a subcomplex whose forms have
well-defined restrictions to boundary strata.  The special differential
forms considered in \S\ref{sectSpecialDifferentialForms} have restrictions
to boundary strata but they are not sufficient to represent
$L^2$-cohomology.  Instead we define in this section a
functor of \emph{quasi-special differential forms},
\begin{equation}
\label{eqnQuasiSpecialDifferentialFormFunctor}
\Rep_\lreg(L_R)\longrightarrow \Complex_\X(\Xhat_R), \qquad E_R \longmapsto
\Adec(\Xhat_R;\EE_R),
\end{equation}
for each $R\in \Pl$ and prove it has the desired properties.  The
definition will be in \S\ref{ssectQuasiSpecialForms}; before that we
define and study quasi-special differential forms $\Adec(\bar{\mathscr
  A}_P^R;\EE_P)$ on $\bar{\mathscr A}_P^R$.

\subsection{Stratification of \boldmath$\bar{\mathscr A}_P^G$}
Recall that $\bar{\mathscr A}_P^G$ is stratified by its $A_P^G$-orbits and
these are indexed by those $Q\in \Pl$ satisfying $P\le Q\le G$.  By
\eqref{eqnCompareScriptAPartialCompactification}, the product decomposition
$\mathscr A_Q^G\times \mathscr A_P^Q$ of the dense stratum $\mathscr A_P^G$
extends to a decomposition $\bar{\mathscr A}_Q^G\times \mathscr A_P^Q$ of
the open star neighborhood of the stratum associated to $Q$; in view of
this we denote the $Q$-stratum $\{o_Q\} \times \mathscr A_P^Q$ or sometimes
simply $\mathscr A_P^Q$.  Define a \emph{special neighborhood} of a
boundary point $z \in \{o_Q\} \times \mathscr A_P^Q$ to be any set of the
form $\bar{\mathscr A}_Q^G(b) \times V^Q$ where $V^Q$ is a relatively
compact neighborhood of $z$ in $\mathscr A_P^Q$.

\subsection{Quasi-special Differential Forms on \boldmath$\bar{\mathscr
    A}_P^R$}%
Let $E_P$ be a locally regular $L_P$-module.  In order to define
$\Adec(\bar{\mathscr A}_P^R;\EE_P)$ we will use induction on $\#\D_P^R$.
If $R=P$ set $\Adec(\bar{\mathscr A}_P^P;\EE_P) = E_P$.  For general $R$,
we may assume by induction that $\Adec(\bar{\mathscr A}_P^Q;\EE_P)$ has
been defined for all $P\le Q <R$.  We will also for simplicity of notation
assume that $R=G$.

Define the sheaf $\Adec(\bar{\mathscr A}_P^G;\EE_P)$ to be the subcomplex
of $\i_{G*}\A(\mathscr A_P^G;\EE_P)$ whose sections over $U \subseteq
\bar{\mathscr A}_P^G$ are those forms $\psi$ on $U\cap \mathscr A_P^G$
which satisfy the following two conditions for all $P \le Q < G$:
\begin{subequations}
\begin{equation}
\label{eqnQuasiSpecialI}
\parbox[t]{.9\textwidth}{%
For every boundary point $z\in U\cap (\{o_Q\} \times \mathscr A_P^Q)$,
there exists a special neighborhood $V= \bar{\mathscr A}_Q^G(b) \times V^Q$
such that $\psi|_{V\cap \mathscr A_P^G} =\sum_i \pr_Q^* \psi_{Q,i} \wedge
(\pr^Q)^*\psi^Q_i$, where $\psi_{Q,i} \in A(\mathscr A_Q^G(b);\EE_P)$ and
$\psi^Q_i\in A(V^Q)$.
}
\end{equation}
Given \eqref{eqnQuasiSpecialI}, the ``restriction'' $\widetilde
m_{QG}(\psi)$ of $\psi$ to $\{o_Q\}\times \mathscr A_P^Q$ may defined
locally as $\sum_i (\psi_{Q,i})_\infty \otimes \psi^Q_i$, where
$(\psi_{Q,i})_\infty\in \i_{G*}\A(\mathscr A_Q^G;\EE_P)_\infty$ denotes the
associated germ.  This is a section of $\i_{Q*}\A(\mathscr
A_P^Q;\i_{G*}\A(\mathscr A_Q^G;\EE_P)_\infty)$.  Now require:
\begin{equation}
\label{eqnQuasiSpecialII}
\parbox[t]{.9\textwidth}{%
$\widetilde m_{QG}(\psi)$ is a section of $\ihat_{Q*}\Adec(\bar{\mathscr
A}_P^Q;\i_{G*}\A(\mathscr A_Q^G;\EE_P)_\infty)$.
}
\end{equation}
\end{subequations}
Condition ~\eqref{eqnQuasiSpecialII} is asserting that $\widetilde
m_{QG}(\psi)$ is a quasi-special form on $\bar{\mathscr A}_P^Q$ (this
notion is well-defined by our inductive hypothesis).

\begin{rem*}
If in \eqref{eqnQuasiSpecialI} we required that $\psi_{Q,i}$ be constant,
we obtain the forms locally lifted from the boundary (see following
\eqref{eqnSpecial}); these are the analogue of special differential forms
in the current context.  If this were satisfied for all $Q$ then
\eqref{eqnQuasiSpecialII} is automatic.  Thus the forms which are locally
lifted from the boundary form a subcomplex of the quasi-special
differential forms on $\bar{\mathscr A}_P^G$.
\end{rem*}

For a given $Q\ge P$, it is not difficult to see that if
~\eqref{eqnQuasiSpecialI} and ~\eqref{eqnQuasiSpecialII} are
satisfied for all $Q'\ge Q$, then we can arrange that
$(\psi_{Q,i})_\infty\in \Adec(\bar{\mathscr A}_Q^G;\EE_P)_\infty$.  We thus
have a morphism
\begin{equation}
\label{eqnRestrictionSplitQuasiSpecialForms}
\widetilde m_{QG}\colon \Adec(\bar{\mathscr A}_P^G;\EE_P) \longrightarrow
\ihat_{Q*}\Adec(\bar{\mathscr A}_P^Q;\Adec(\bar{\mathscr
A}_Q^G;\EE_P)_\infty)
\end{equation}
such that $\i_Q^*(\widetilde m_{QG})$ is an isomorphism.  Note that the
analogue of \eqref{eqnComposeRestrictionSpecialForms} does not hold.

If $E_P$ is a regular $L_P$-module, define
$\widetilde{\A}_{(2),\textup{sp}}(\bar{\mathscr A}_P^G;\EE_P,\wgt_P) =
\Ltwo(\bar{\mathscr A}_P^G;\EE_P,\wgt_P)\cap \Adec(\bar{\mathscr
A}_P^G;\EE_P)$.  The mapping $\widetilde m_{QG}$ of
\eqref{eqnRestrictionSplitQuasiSpecialForms} restricts to
\begin{equation}
\label{eqnRestrictionSplitLTwoQuasiSpecialForms}
\widetilde m_{QG}\colon \widetilde{\A}_{(2),\textup{sp}}(\bar{\mathscr
  A}_P^G;\EE_P,\wgt_P) \longrightarrow 
\ihat_{Q*}\Adec(\bar{\mathscr A}_P^Q;
  \widetilde{\A}_{(2),\textup{sp}}(\bar{\mathscr A}_Q^G;\EE_P,\wgt_Q)_\infty)
\end{equation}
since $h_P|_{\mathscr A_Q^G\times V^Q} \sim \pr_Q^* h_Q$ and again
$\i_Q^*(\widetilde m_{QG})$ is an isomorphism.  For general $R\neq G$,
$\widetilde m_{QR}$ and its restriction to $L^2$-forms are defined
analogously to \eqref{eqnRestrictionSplitQuasiSpecialForms} and
\eqref{eqnRestrictionSplitLTwoQuasiSpecialForms}.

\begin{prop}
\label{propSplitQuasiSpecialQuasiIsomorphism}
Let $E_P$ be a locally regular $L_P$-module.
\begin{enumerate}
\item
\label{itemSplitQuasiSpecialQuasiIsomorphism}
The inclusion map is a quasi-isomorphism
\begin{equation*}
\Adec(\bar{\mathscr A}_P^G;\EE_P)_\infty \longrightarrow \i_{G*}\A(\mathscr
A_P^G;\EE_P)_\infty.
\end{equation*}
\item
\label{itemSplitQuasiSpecialLtwoQuasiIsomorphism}
If $E_P$ is a regular $L_P$-module and $\EE_P$ is given the
fiber metric coming from an admissible inner product, the inclusion map is
a quasi-isomorphism
\begin{equation*}
\widetilde{\A}_{(2),\textup{sp}}(\bar{\mathscr A}_P^G;\EE_P,\wgt_P)_\infty
\longrightarrow \Ltwo(\bar{\mathscr A}_P^G;\EE_P,\wgt_P)_\infty.
\end{equation*}
\end{enumerate}
\end{prop}
\begin{rem*}
If $E_P$ is irreducible and $\xi_{E_P}+\hsr_P$ as in
\eqref{eqnQQtildeDefinition} is either regular or non-dominant with respect
to $\D_P$, the local $L^2$-cohomology in part
\itemref{itemSplitQuasiSpecialLtwoQuasiIsomorphism} is finite dimensional
and is represented by constant forms \cite[(4.51)]{refnZuckerWarped}.  We
need to consider the general situation however where the local
$L^2$-cohomology may be infinite dimensional.
\end{rem*}

The proof will appear below in
\S\ref{ssectProofpropSplitQuasiSpecialQuasiIsomorphism} after some
preliminaries on homotopy operators.  In order to use coordinates defined
by roots, we assume a basepoint has been chosen and identify $\bar{\mathscr
A}_P^G$ with $\Abar_P^G$ for the remainder of the section.

\subsection{Basic Homotopy Operator}
\label{ssectBasicHomotopyOperator}
{}From now through \S\ref{ssectLTwoHomotopy} we consider a fixed $Q\ge P$.
We would like to have a homotopy formula between a form $\psi$ on $A_P^G$
and a form that satisfies \eqref{eqnQuasiSpecialI} and
\eqref{eqnQuasiSpecialII} near the $Q$-stratum.  We write $a=a_Q a^Q\in
A_P^G$ according to the decomposition $A_P^G = A_Q^G \times A_P^Q$.

Fix $c\in A_P^Q$ and let $\pi_c\colon A_P^G\to A_Q^G\times \{c\} \subseteq
A_P^G$ be the projection $a_Q a^Q\mapsto a_Q c$.  Order $\D_P^Q = \{
\al_1,\dots,\al_n\}$ and let $x^i(a) = \log a^{\al_i}$ be the corresponding
coordinates on the $A_P^Q$ factor.  Let $c^i = x^i(c)$ for $1\le i \le n$
and set $(A_P^G)_{c,i} = \{\, a\in A_P^G\mid x^j(a) = c^j \text{ for }1\le
j\le i\,\}$.  Let $\pi_{c,i}\colon A_P^G \to (A_P^G)_{c,i} \subseteq A_P^G$
be the coordinate projection; note that $\pi_{c,0}=\id_{A_P^G}$ and
$\pi_{c,n}=\pi_c$.  As in \cite{refnCheeger} we have homotopy and
projection operators
\begin{equation}
\label{eqnQHomotopy}
H_{Q,c}\psi = \sum_{i=1}^n \int_{c^i}^{x^i}
\iota_{\frac{\partial}{\partial x^i}} \pi_{c,i-1}^*\psi \qquad \text{and}
\qquad P_{Q,c}\psi = \pi_c^*\psi,
\end{equation}
where the integral is integrating with respect to the
$i^\text{th}$-coordinate.  One calculates that
\begin{equation}
\label{eqnQHomotopyFormula}
dH_{Q,c}\psi + H_{Q,c}d\psi = \psi - P_{Q,c}\psi.
\end{equation}

\subsection{Homotopy of Forms on Neighborhoods}
\label{ssectTruncatedNeighborhoods}
We wish to apply $H_{Q,c}$ to forms on $U\cap A_P^G$, where $U$ belongs to
a fundamental system of neighborhoods of $o_P\in \Abar_P^G$.  We will use
the standard fundamental system of neighborhoods given by
\begin{equation}
\label{eqnShiftedCone}
\Abar_P^G(s) = \{\,a\in \Abar_P^G \mid a^\al>s \text{ for all $\al \in
  \D_P$}\,\}
\end{equation}
for $s\ge 1$; we set $A_P^G(s) = \Abar_P^G(s)\cap A_P^G$.  The difficulty
is that no matter what $c\in A_P^Q(s)$ is chosen (unless $\D_P^Q$ is
orthogonal to $\D_P\setminus \D_P^Q$) there will exist points $a\in
A_P^G(s)$ such that $\pi_c(a) \notin A_P^G(s)$.  Thus if $\psi$ is a form
on $A_P^G(s)$, equation \eqref{eqnQHomotopy} does not define $H_{Q,c}\psi$
and $P_{Q,c}\psi$ as forms on $A_P^G(s)$.  Nonetheless we have the
following lemma:

\begin{lem}
\label{lemHomotopyFormulaNeighborhoods}
For all $s \ge 1$ and for all $c\in A_P^Q(s)$, equation
\eqref{eqnQHomotopy} defines operators
\begin{equation*}
H_{Q,c}\colon
A(A_P^G(s); \EE_P) \to A(A_P^G(m_cs); \EE_P)[-1]
\end{equation*}
and
\begin{equation*}
P_{Q,c}\colon
A(A_P^G(s); \EE_P) \to A(A_P^G(m_cs); \EE_P),
\end{equation*}
where $m_c = \max_{\al\in \D_P^Q} c^\alpha/s > 1$, and the homotopy formula
\eqref{eqnQHomotopyFormula} holds.
\end{lem}
\begin{proof}
  We need to verify that $H_{Q,c}\psi$ and $P_{Q,c}\psi$ are defined at
  $a=a_Q a^Q\in A_P^G(m_cs)$; note that in this case $c^\al < (a^Q)^\al$
  for all $\al\in\D_P^Q$.  It suffices to check that if $b\in A_P^Q$ is any
  element such that $c^\al \le b^\al \le (a^Q)^\al$ for all $\al \in
  \D_P^Q$, then $a_Qb\in A_P^G(s)$.  For $\al\in \D_P^Q$, $(a_Qb)^\al =
  b^\al \ge c^\al > s$.  On the other hand, for $\d \in \D_P \setminus
  \D_P^Q$, write $\d|_{A_P^Q} = \sum_{\al\in\D_P^Q} \langle \d,
  \b_\al^Q{}\spcheck \rangle \al$ where $\langle \d, \b_\al^Q{}\spcheck
  \rangle \le 0$.  Then $(a_Qb)^\d > m_c s (b/a^Q)^{\d} = m_c
  s\prod_{\al\in\D_P^Q} (b/a^Q)^{\langle \d, \b_\al^Q{}\spcheck \rangle \al}
  \ge m_c s > s$.
\end{proof}

\subsection{Cutoff Homotopy Operator}
If a form $\psi$ on $A_P^G(s)$ already satisfies \eqref{eqnQuasiSpecialI}
and \eqref{eqnQuasiSpecialII} near the $R$-stratum with $R\neq Q$, the same
may not be true for $H_{Q,c}\psi$ and $P_{Q,c}\psi$.  We will remedy this
by multiplying $H_{Q,c}$ by a certain cutoff function in order to restrict
its effect to points near the $Q$-stratum.  For later use in
\S\ref{ssectLTwoHomotopy} we need to allow some flexibility in the cutoff
function.

Let us call a smooth function $g\colon \RR \to [1,\infty)$
\emph{$\epsilon$-admissible} if
\begin{enumerate}
\item $g(r)$ is monotonically increasing to $\infty$, and
\item $g'(r)< \epsilon$ for all $r\in \RR$.
\end{enumerate}
Clearly $\epsilon$-admissible functions exist for any $\epsilon>0$.  For an
$\epsilon$-admissible function $g$, define
\begin{equation}
\label{eqnChiDefinition}
\chi_g= \chi_g(a) = \cutoff\bigl(\max_{\al\in\D_P^Q} \log{(a^Q)^\al} - g(
  \min_{\g\in\D_Q} \log a_Q^{\g})\bigr)
\end{equation}
where $\cutoff(x)$ is a smooth cutoff function such that $\cutoff(x)=1$ for
$x\le -1$ and $\cutoff(x)=0$ for $x\ge 0$.  Note that $\chi_g$ is not
smooth; one can rectify this by replacing $\max$ and $\min$ by appropriate
smooth approximations.  (Alternatively one could use piecewise smooth forms
and the distribution exterior derivative throughout.)

Define a cutoff homotopy operator and projection,
\begin{gather}
\label{eqnCutoffHomotopyOperator}
H_{Q,c,g}\psi = \chi_g H_{Q,c}\psi,\\
\label{eqnProjectionDefinition}
P_{Q,c,g}\psi = (1- \chi_g) \psi + \chi_g P_{Q,c}\psi + d\chi_g\wedge
H_{Q,c}\psi,
\end{gather}
where $H_{Q,c}$ and $P_{Q,c}$ are obtained from Lemma
~\ref{lemHomotopyFormulaNeighborhoods}.  It is straightforward to calculate
from \eqref{eqnQHomotopyFormula} that for $\psi \in A(A_P^G(s); \EE_P)$,
\begin{equation}
\label{eqnGHomotopyFormula}
dH_{Q,c,g}\psi + H_{Q,c,g} d\psi = \psi|_{A_P^G(m_cs)} - P_{Q,c,g}\psi.
\end{equation}

\begin{lem}
\label{lemHomotopyWithPreservation}
Let $c\in A_P^Q(s)$ and let $g$ be any $\epsilon$-admissible function.  For
any $\psi \in A(A_P^G(s); \EE_P)$, the form $P_{Q,c,g} \psi$ satisfies
\eqref{eqnQuasiSpecialI} and \eqref{eqnQuasiSpecialII} near the boundary
component corresponding to $Q$.  If $\epsilon>0$ is sufficiently small
\textup(depending only on $G$\textup) and if $\psi$ already satisfies
\eqref{eqnQuasiSpecialI} and \eqref{eqnQuasiSpecialII} near the boundary
component corresponding to some $R\nless Q$, then so do $H_{Q,c,g}\psi$ and
$P_{Q,c,g}\psi$.
\end{lem}

\begin{proof}
In a special neighborhood $V= \Abar_Q^G(b) \times V^Q$ of a point on the
$Q$-stratum, $\log a^{Q\al}$ is bounded for all $\al\in\D_P^Q$ and by
making $b$ sufficiently regular we can arrange that $\log a_Q^\g$ is
arbitrarily large for all $\g\in \D_Q$.  Since $g$ tends to infinity, we
can arrange that $\chi_g|_V \equiv 1$ and thus $P_{Q,c,g}\psi|_V = \pi_c^*\psi$.
This proves \eqref{eqnQuasiSpecialI} for $P_{Q,c,g}\psi$ near the $Q$-stratum.
Furthermore this shows that the restriction $\widetilde m_{QG}(P_{Q,c,g}\psi)$ is
the function on $\{o_Q\}\times A_P^Q$ which has the constant value
$(\psi|_{A_P^G(s)\cap(A_Q^G\times\{c\})})_\infty \in \i_{G*}\A(\mathscr
A_Q^G;\EE_P)_\infty$.  Thus \eqref{eqnQuasiSpecialII} holds as well.

Now consider $R\nless Q$ such that \eqref{eqnQuasiSpecialI} and
\eqref{eqnQuasiSpecialII} hold near the $R$-stratum.  We will prove the
final assertion by induction on $\#\D_P^G$.  We will just treat $H_{Q,c,g}\psi$
since the argument for $P_{Q,c,g}\psi$ is identical.  Let $V = \Abar_R^G(b)
\times V^R$ be a special neighborhood of a point on the $R$-stratum and
write, as in \eqref{eqnQuasiSpecialI},
\begin{equation*}
\psi|_{V\cap \mathscr A_P^G} =\sum_i \pr_R^*\psi_{R,i}\wedge
(\pr^R)^*\psi^R_i.
\end{equation*}

First assume $R\ge Q$.  Then we can assume that $V^R$ decomposes into
$V_Q^R \times V^Q$ according to $A_P^R = A_Q^R \times A_P^Q$ and that both
factors are relatively compact.  By enlarging $V^Q$ we may also assume that
$c\in V^Q$.  For $a\in V\cap A_P^G$, decompose $a = a_R a_Q^R a^Q$
according to these decompositions; note that $a^R = a_Q^R a^Q$ and $a_Q =
a_R a_Q^R$.  Since $V_Q^R$ is relatively compact, $(a_Q^R)^\g$ belongs to a
relatively compact subset of $(0,\infty)$ for all $\g\in\D_Q$.  Thus if $\g
\in \D_Q^R$ then $a_Q^\g = (a_Q^R)^\g$ is bounded.  If $\g \in
\D_Q\setminus \D_Q^R$ then $a_Q^\g = a_R^\g (a_Q^R)^\g$ can be made
arbitrarily large by making $b$ more regular (since $a_R \in A_R^G(b)$).
Thus $\min_{\g\in\D_Q} \log a_Q^{\g} = \min_{\g\in\D_Q^R} \log
(a_Q^R)^{\g}$ and hence $\chi_g(a)$ depends only on $a^R$.  This means that
\begin{equation*}
(H_{Q,c,g}\psi)|_{V\cap \mathscr A_P^G} =\sum_i \pr_R^*\psi_{R,i}\wedge
(\pr^R)^*(H_{Q,c,g}\psi^R_i).
\end{equation*}
Thus \eqref{eqnQuasiSpecialI} holds for $H_{Q,c,g}\psi$ near the
$R$-stratum and also $\widetilde m_{RG}(H_{Q,c,g}\psi) = \linebreak H_{Q,c,g}
\widetilde m_{RG}(\psi)$.  But the lemma applies by induction to
$\widetilde m_{RG}(\psi)$, which is quasi-special, and implies that
$H_{Q,c,g} \widetilde m_{RG}(\psi)$ is quasi-special.  Thus
\eqref{eqnQuasiSpecialII} holds.

Now assume $R\nleq Q$ and $R\ngeq Q$.  The first of these conditions
implies that there exists $\al' \in \D_P^R$ such that $\al' \notin \D_P^Q$;
set $\g = \al'|_{A_Q}$.  For $a\in V\cap A_P^G$ we have $a_Q^\g = a^{\al'}
(a^Q)^{-\al'} = (a^R)^{\al'} (a^Q)^{-\al'}$.  But since $a^R \in V^R$ and
$V^R$ is relatively compact, we may estimate that
\begin{equation*} 
  \begin{split}
    \log a_Q^\g  &\le C+ \langle -\al', \tau_P^Q{}\spcheck\rangle
    \max_{\al\in\D_P^Q} \log (a^Q)^\al  \\
    &\le C + M \cdot \max_{\al\in\D_P^Q}\log (a^Q)^\al,
  \end{split}
\end{equation*}
where $C>0$ and $M = \max_{\d\in\D_P\setminus \D_P^Q} \langle -\d,
\tau_P^Q{}\spcheck\rangle \ge 0$.    Thus for $a\in V$, 
\begin{equation}
\label{eqnGEstimate}
g(\min_{\g\in\D_Q} \log a_Q^{\g}) \le
g(C) + \epsilon M \cdot
\max_{\al\in\D_P^Q}\log (a^Q)^\al.
\end{equation}
Next, the condition $R\ngeq Q $ implies that there exists $\al''\in \D_P^Q$
such that $\al''\notin \D_P^R$ and hence  $(a^Q)^{\al''} = a^{\al''} =
a_R^{\al''} (a^R)^{\al''}$.  For $a\in V$ we have $a_R \in \Abar_R^G(b)$
and $a^R \in V^R$; since $V^R$ is relatively compact, this implies that
$(a^Q)^{\al''}$ can be made arbitrarily large by making $b$ more regular.
In particular we can arrange that
\begin{equation*}
\label{eqnLogAQBound}
g(C) \le \frac 12 \max_{\al\in\D_P^Q}\log (a^Q)^\al, \qquad \text{for all
  $a\in V$}.
\end{equation*}
This bounds the first term of \eqref{eqnGEstimate}, while we can choose
$\epsilon>0$ to arrange that the second term of \eqref{eqnGEstimate} is $\le
\frac12 \max_{\al\in\D_P^Q}\log (a^Q)^\al$ as well.  Thus $\chi_g|_V \equiv
0$.  Therefore $H_{Q,c,g}\psi|_V = 0$ and $P_{Q,c,g}\psi|_V = \psi$ so
\eqref{eqnQuasiSpecialI} and \eqref{eqnQuasiSpecialII} trivially hold.
\end{proof}

\subsection{Homotopy of \boldmath$L^2$ Forms}
\label{ssectLTwoHomotopy}
Now assume that $E_P$ is regular with an admissible inner product.  In
general the homotopy operator $H_{Q,c,g}$ for any fixed choice of $g$ will not be
bounded on $L^2$.  We will get around this by choosing $g$ depending on
each given $\psi \in A_{(2)}(A_P^G(s); \EE_P,\wgt_P)$.

Set $A_P^G(s)_{c,i} = A_P^G(s) \cap (A_P^G)_{c,i} = \{\, a\in
A_P^G(s) \mid x^j(a) = c^j \text{ for }1\le j\le i\,\}$ and
\begin{equation}
\sigma_Q = \langle \tau_Q^G,\tau_Q^G{}\spcheck\rangle^{-1} \tau_Q^G.
\end{equation}
For an $\epsilon$-admissible function $g$ we will also denote by $g$ the
weight function $a\mapsto g(\log a_Q^{\sigma_Q})$ on $A_P^G(s)$.

\begin{lem}
\label{lemGFunction}
Let $\{\psi_\mu\}$ be a sequence in $A(A_P^G(s)_{c,i}; \EE_P)$, where $c\in
A_P^Q(s)$ and $0\le i \le n$.  Assume that $\|\psi_\mu\|_{\wgt_P} < \infty$
for all $\mu$ and that $\{\psi_\mu\}$ is convergent in this norm.  For any
$\epsilon>0$ there exists an $\epsilon$-admissible function $g$ such that
$\|\psi_\mu\|_{g\wgt_P} < \infty$ for all $\mu$ and, after passing to a
subsequence, $\{\psi_\mu\}$ is convergent in this new norm.
\end{lem}

\begin{proof}
We assume that $i=0$, that is, $\{\psi_\mu\}$ is a sequence in
$A(A_P^G(s); \EE_P)$; the general case is identical.  Assume
the sequence starts with $\psi_0=0$ and let $\psi = \lim_{\mu\to\infty}
\psi_\mu$.  If we set $A_P^G(s)^t = \{\, a\in A_P^G(s) \mid \log
a_Q^{\sigma_Q}=t\,\}$, then by Fubini's theorem
\begin{equation}
\label{eqnPhiLtwoNorm}
\|\psi_\mu-\psi\|_{\wgt_P}^2 = \int_{0}^\infty \biggl(
\int_{A_P^G(s)^t} |\psi_\mu-\psi|^2 \wgt_P(a)^2\, da \biggr) \, dt
\end{equation}
where $da$ is the induced measure on $A_P^G(s)^t$.  Let $f_{\mu}(t)$
denote the inner integral in \eqref{eqnPhiLtwoNorm}.  Note that
$\lim_{\mu\to \infty} \int_t^\infty f_{\mu}\,dt = 0$ for all $t$ and that
$\lim_{t\to \infty} \int_t^\infty f_{\mu}\,dt = 0$ for all $\mu$.  Thus we
can find a sequence $t_1< t_2< \dots< t_k\to\infty$ such that
$\int_{t_k}^\infty f_{\mu}\,dt \le (k+1)^{-3}$ for all $\mu$, $k$ and
$t_{k+1} \ge t_k + (2/\epsilon)$ for all $k$.  (Given $k$, the inequality
is automatic for large $\mu$ and then $t_k$ can be made sufficiently large
to accommodate the rest.)  It follows that $\sum_{k=1}^\infty (k+1)
\int_{t_k}^{t_{k+1}} f_\mu\,dt < \infty$ for all $\mu$.  Let $g(t)$ be a
smooth monotonic function such that $g(t)=1$ for $t\le t_1$, $g(t_k) =
\sqrt{k}$ for $k\ge 1$, and $g'(t)\le \epsilon$.  Then
$\|\psi_\mu-\psi\|_{g\wgt_P} < \infty$ for all $\mu$ which implies
$\|\psi\|_{g\wgt_P} < \infty$ (set $\mu=0$) and hence
$\|\psi_\mu\|_{g\wgt_P} < \infty$ for all $\mu$.

For the final assertion, choose a subsequence such that $\int_0^\infty
f_{\mu}\,dt \le \mu^{-2}$ for $\mu\ge 1$ and hence $\lim_{\mu\to \infty}
\mu \int_t^\infty f_{\mu}\,dt = 0$ for all $t$.  We can then choose $t_k$
in the argument above such that $\|\psi_\mu-\psi\|_{g\wgt_P} < \mu^{-1}$
for $\mu\ge 1$.
\end{proof}

\begin{rem}
\label{remGFunction}
The lemma continues to hold (with the same proof) if $\psi$ or $\psi_\mu$
are measurable, not necessarily smooth, $L^2$ forms.  Furthermore, given a
finite number of convergent sequences, a single function $g$ can be found
so that the conclusion of the lemma will hold for all the sequences.  Our
main interest will be when these sequences are actually constant.
\end{rem}

\begin{lem}
\label{lemExistenceOfAdmissibleFunction}
Let $\psi \in A_{(2)}(A_P^G(s); \EE_P,\wgt_P)$.  For any $\epsilon>0$ and
for almost every $c\in A_P^Q(s)$, there exists an $\epsilon$-admissible
function $g$ such that $H_{Q,c,g}\psi$, $P_{Q,c,g}\psi\in
A_{(2)}(A_P^G(m_cs); \EE_P,\wgt_P)$.
\end{lem}

\begin{proof}
Pick $c\in A_P^Q(s)$ such that for all $0\le i \le n$ both
$\psi|_{A_P^G(s)_{c,i}}$ and $d\psi|_{A_P^G(s)_{c,i}}$ are $L^2$ with
weight $\wgt_P$; this condition is satisfied by almost every $c$.  Apply
Lemma ~\ref{lemGFunction} (see also Remark ~\ref{remGFunction}) to find an
$\epsilon$-admissible function $g$ such that all the above forms are $L^2$
with weight $g\wgt_P$.  For $a\in \supp \chi_g$ we have the estimate
\begin{equation}
x^i(a) = \log a ^{\al_i} \le \max_{\al\in\D_P^Q} \log a^{Q\al} \le g(
  \min_{\g\in\D_Q} \log a_Q^{\g}) \le
  g(\log a_Q^{\sigma_Q})
\end{equation}
where the last inequality comes from the estimate
\begin{equation*}
a_Q^{\tau_Q^G} =  a_Q^{\sum_{\g\in\D_Q}\langle\tau_Q^G,\b_\g\spcheck\rangle
  \g}
\ge (\min_{\g\in\D_Q} a_Q^{\g})^{\sum_\g
  \langle\tau_Q^G,\b_\g\spcheck\rangle} 
=  (\min_{\g\in\D_Q} a_Q^{\g})^{\langle\tau_Q^G,\tau_Q^G{}\spcheck\rangle}.
\end{equation*}
On the other hand, for any $a \in A_P^G(s)$ we have the estimate
\begin{equation}
c^i \lesssim g( \log s^{1+\min_{\d} \langle -\d,
\tau_P^Q{}\spcheck\rangle} ) \le  g( \min_{\g\in\D_Q} \log
  a_Q^{\g}) \le  g(\log a_Q^{\sigma_Q}).
\end{equation}
(Here and below we will use the notation $p \lesssim q$ to indicate $p \le
C q$ where $C>0$ is a sufficiently large constant depending only on $P$,
$Q$, $c$, $g$, and $s$.)  We also observe that $h_P(a)^2 = h_P(a_Q)^2
\prod_i e^{-d_ix^i}$ where $d_i>0$.  From these facts we can estimate
\begin{equation}
\label{eqnHomotopyBounded}
  \begin{split}
    \|H_{Q,c,g}\psi\|^2_{\wgt_P} &\lesssim \sum_{i=1}^n \int_{A_P^G(m_cs)}
    \biggl(\chi_g \int_{c^i}^{x^i} |\pi_{c,i-1}^*\psi|\biggr)^2
    \wgt_P(a)^2 \,dV \\
    &\lesssim \sum_{i=1}^n \int_{A_P^G(m_cs) \cap \supp \chi_g}
    \biggl(\int_{c^i}^{x^i}
    |\pi_{c,i-1}^*\psi|^2\biggr) g(\log a_Q^{\sigma_Q}) \wgt_P(a)^2 \,dV \\
    &\lesssim \sum_{i=1}^n \int_{A_P^G(s) \cap \supp \chi_g}
    |\pi_{c,i-1}^*\psi|^2 g(\log a_Q^{\sigma_Q}) \wgt_P(a)^2
    \,dV \\ 
    &\lesssim \sum_{i=1}^n  \| \psi|_{A_P^G(s)_{c,i-1}}
    \|_{g\wgt_P}^2
    < \infty
  \end{split}
\end{equation}
where the inner integrals are with respect to the $i^\text{th}$-coordinate
and for the third line we use Fubini's theorem.  We may similarly prove
$\|H_{Q,c,g}d\psi\|_{\wgt_P}$, $\|\chi_gP_{Q,c}\psi \|_{\wgt_P}$, and $\|
d\chi_g \wedge H_{Q,c}\psi\|_{\wgt_P}$ are finite which by
\eqref{eqnProjectionDefinition} shows that
$\|P_{Q,c,g}\psi\|_{\wgt_P}<\infty$ as well.  Finally
\eqref{eqnGHomotopyFormula} shows that $\|dH_{Q,c,g}\psi\|_{\wgt_P} <
\infty$.  The same argument applied to $d\psi$ shows that
$\|dH_{Q,c,g}d\psi\|_{\wgt_P} < \infty$ and hence by
\eqref{eqnGHomotopyFormula} that $\|dP_{Q,c,g}\psi\|_{\wgt_P} < \infty$.
\end{proof}
\begin{rem}
Suppose that $\psi_\mu \to \psi$ in the graph norm on $A_{(2)}(A_P^G(s);
\EE_P,\wgt_P)$.  Then for all $i$ and almost every $c\in A_P^Q(s)$,
$\psi_\mu |_{A_P^G(s)_{c,i}} \to \psi|_{A_P^G(s)_{c,i}}$ in the graph norm.
Hence the argument of the lemma shows that for any $\epsilon>0$ and for
almost every $c\in A_P^Q(s)$, there exists an $\epsilon$-admissible
function $g$ such that after passing to a subsequence, $H_{Q,c,g}\psi_\mu
\to H_{Q,c,g}\psi$ and $P_{Q,c,g}\psi_\mu \to P_{Q,c,g}\psi$ in the graph
norm of $A_{(2)}(A_P^G(m_cs); \EE_P,\wgt_P)$.  Again the result continues to
hold for measurable forms which together with their exterior derivatives in
the distribution sense are $L^2$.
\end{rem}

\subsection{Proof of Proposition
  ~\ref{propSplitQuasiSpecialQuasiIsomorphism}}
\label{ssectProofpropSplitQuasiSpecialQuasiIsomorphism}
Choose a total ordering $Q_0=P$, $Q_1$, $Q_2$, \dots, $Q_N$, $Q_{N+1}=G$ of
the parabolic $\QQ$-subgroups containing $P$ which is compatible with the
partial order; such a total ordering exists.  Let
$\widetilde{\A}_{(2),\textup{sp},i}(\Abar_P^G;\EE_P,\wgt_P)$ be the
subcomplex of $ \Ltwo(\bar{\mathscr A}_P^G;\EE_P,\wgt_P)$ whose sections
satisfy \eqref{eqnQuasiSpecialI} and \eqref{eqnQuasiSpecialII} near points
of the $Q_j$-strata for all $j \ge i$, and let
$\Atilde_{(2),\textup{sp},i}(U;\EE_P,\wgt_P)$ denote its sections over
$U\subseteq \Abar_P^G$.

For part ~\itemref{itemSplitQuasiSpecialLtwoQuasiIsomorphism} we will
show for any $0\le i \le N$ that the map of local cohomology
\begin{equation}
\label{eqnQuasiSpecialLocalCohomology}
H(\widetilde{\A}_{(2),\textup{sp},i}(\Abar_P^G;\EE_P)_\infty) \to
H(\widetilde{\A}_{(2),\textup{sp},i+1}(\Abar_P^G;\EE_P)_\infty)
\end{equation}
is an isomorphism.  For all $s\ge 1$, set $C_{(2),i}(s) =
\Atilde_{(2),\textup{sp},i}(\Abar_P^G(s);\EE_P,\wgt_P)$ and consider
\begin{equation*}
\entrymodifiers={+!!<0pt,\fontdimen22\textfont2>}
\xymatrix {
{C_{(2),i}(s)} \ar@{^{ (}->}[r] \ar[d] & {C_{(2),i+1}(s)}
\ar[d] \\
{C_{(2),i}(2s)} \ar@{^{ (}->}[r] & {C_{(2),i+1}(2s)}
}
\end{equation*}
where the vertical maps are induced by restriction.  If $[\psi]\in
H(C_{(2),i+1}(s))$, use Lemmas ~\ref{lemHomotopyWithPreservation} and
\ref{lemExistenceOfAdmissibleFunction} (applied with $Q=Q_i$ and $R$
ranging over $Q_j$ for all $j>i$) and \eqref{eqnGHomotopyFormula} to find
$c\in A_P^Q(s) \cap A_P^Q(1/(2s))^{-1}$ and $g$ such that
$dH_{Q_i,c,g}\psi = \psi|_{A_P^G(2s)} - P_{Q_i,c,g}\psi$ with
$P_{Q_i,c,g}\psi \in C_{(2),i}(2s)$ and $H_{Q_i,c,g}\psi \in
C_{(2),i+1}(2s)$.  This shows surjectivity in
\eqref{eqnQuasiSpecialLocalCohomology}.  Similarly, if $[\eta]\in
H(C_{(2),i}(s))$ and $\eta=d\psi$ with $\psi \in C_{(2),i+1}(s)$, we find
that $\eta|_{A_P^G(2s)} = dH_{Q_i,c,g}\eta + dP_{Q_i,c,g}\psi$ with
$H_{Q_i,c,g}\eta$, $P_{Q_i,c,g}\psi \in C_{(2),i}(2s)$.  This shows
$[\eta]=0$ and hence injectivity.

The proof of part ~\itemref{itemSplitQuasiSpecialQuasiIsomorphism} is
simpler since any $c$ and any $\epsilon$-admissible $g$ will suffice.  \qed

\subsection{Quasi-special Differential Forms on \boldmath$\Xhat$}
\label{ssectQuasiSpecialForms}
We now define $\Adec(\Xhat;\EE)$, the quasi-special differential forms on
$\Xhat$.  For a locally regular representation $E$ of $G$, a section of
$\Adec(\Xhat;\EE)$ over $U\subseteq \Xhat$ is an element $\eta \in A(U\cap
X;\EE)$ satisfying the following condition for all $Q\le G$:
\begin{equation}
\parbox[t]{.9\textwidth}{%
For every boundary point $x\in U\cap X_Q$, there exists a special
neighborhood $V= \nil(\bar{\mathscr A}_Q^G(b) \times O_Q) \subseteq U$ of
$x$ (see \eqref{eqnSpecialNeighborhood}) such that $\eta|_{V\cap X}$ is a
sum of terms $\pr_Q^* \psi\wedge (\pr^Q)^*\o$, where $\psi\in
\Atilde\sp(\bar{\mathscr A}_Q^G(b);\EE)$ and $\o \in A(O_Q)$ is
$N_Q(\RR)$-invariant.
}
\end{equation}
A special differential form as in \S\ref{ssectRealization} satisfies this
condition with $\psi$ constant.  Thus $\Asp(\Xhat;\EE)$ is a subcomplex of
$\Adec(\Xhat;\EE)$.  Set $\widetilde{\A}_{(2),\textup{sp}}(\Xhat;\EE) =
\Ltwo(\Xhat;\EE)\cap \Adec(\Xhat;\EE)$.
\begin{prop}
\label{propQuasiSpecialQuasiIsomorphism}
Let $E$ be a locally regular $G$-module.
\begin{enumerate}
\item
\label{itemQuasiSpecialQuasiIsomorphism}
The inclusion map is a quasi-isomorphism
\begin{equation*}
\Adec(\Xhat;\EE) \longtildearrow \i_{G*}\A(X;\EE).
\end{equation*}
\item
\label{itemQuasiSpecialLtwoQuasiIsomorphism}
If $E$ is a regular $G$-module and $\EE$ is given the fiber metric coming
from an admissible inner product, then the inclusion map is a
quasi-isomorphism
\begin{equation*}
\widetilde{\A}_{(2),\textup{sp}}(\Xhat;\EE) \longtildearrow
\Ltwo(\Xhat;\EE).
\end{equation*}
\end{enumerate}
\end{prop}

\begin{proof}
Consider the diagram of maps between stalks at a point $x\in X_P$:
\begin{equation*}
\xymatrix {
{\A_{(2)}(\bar{\mathscr A}_P^G; \HH(\n_P;E), \wgt_P)_\infty}
\ar@{^{ (}->}[r]  & {\Ltwo(\Xhat;\EE)_x} \\
{\widetilde{\A}_{(2),\textup{sp}}(\bar{\mathscr A}_P^G; \HH(\n_P;E), \wgt_P
  )_\infty} \ar@{^{ (}->}[r]  
\ar@{^{ (}->}[u] &
{\widetilde{\A}_{(2),\textup{sp}}(\Xhat;\EE)_x.} \ar@{^{ (}->}[u]
}
\end{equation*}
Zucker \cite[(4.24)]{refnZuckerWarped} defines the top inclusion and shows
that it is a quasi-isomorphism.  The proof involves constructing a
projection mapping in the opposite direction (the composition of harmonic
projection in the $O_P$ factor of a special neighborhood, projection to
$N_P(\RR)$-invariant forms, and harmonic projection in the complex
$\bigwedge \n_P^*\otimes E$) and a bounded homotopy operator between it and
the identity.  These operators all preserve the condition that a form is
quasi-special, so the bottom inclusion is also a quasi-isomorphism.  The
left inclusion is a quasi-isomorphism by Proposition
~\ref{propSplitQuasiSpecialQuasiIsomorphism}%
\itemref{itemSplitQuasiSpecialLtwoQuasiIsomorphism} and hence the right
inclusion is as well.  This proves part
\itemref{itemQuasiSpecialLtwoQuasiIsomorphism}; the proof of part
\itemref{itemQuasiSpecialQuasiIsomorphism} is similar using Proposition
~\ref{propSplitQuasiSpecialQuasiIsomorphism}%
\itemref{itemSplitQuasiSpecialQuasiIsomorphism}.
\end{proof}

\subsection{The Quasi-special Differential Forms Functor}
If we apply the construction of \S\ref{ssectQuasiSpecialForms} to each
$\Xhat_R$, we obtain the desired functors
\eqref{eqnQuasiSpecialDifferentialFormFunctor}.

\begin{prop}
\label{propRestrictionTildeSpecialForms}
Let $E_Q$ be a locally regular $L_Q$-module.
\begin{enumerate}
\item
\label{itemTildeSpecialQuasiIsomorphism}
The inclusion morphism yields a natural quasi-isomorphism
\begin{equation*}
\Asp(\Xhat_Q;\EE_Q) \longtildearrow \Adec(\Xhat_Q;\EE_Q).
\end{equation*}

\item
\label{itemRestrictionTildeSpecialForms}
For $P\le Q \in \Pl$, the ``restriction'' morphism $k_{PQ}$ of special
differential forms \textup(see \eqref{eqnRestrictionSpecialForms}\textup)
extends to a natural morphism
\begin{equation*}
\widetilde k_{PQ}:\Adec(\Xhat_Q;\EE_Q) \longrightarrow
\ihat_{P*}\Adec(\Xhat_P; \Adec(\bar{\mathscr A}_P^Q; \HH(\n_P^Q; E_Q))_\infty)
\end{equation*}
such that $\ihat_P^* (\widetilde k_{PQ})$ is a quasi-isomorphism.  This
morphism satisfies
\begin{equation}
\label{eqnComposeRestrictionQuasiSpecialForms}
\ihat_{Q*}(\widetilde k_{PQ})\circ \widetilde k_{QR} = \ihat_{P*}\Adec(\Xhat_P;
  \widetilde m_{QR}) \circ \widetilde k_{PR}.
\end{equation}

\item
\label{itemRestrictionLTwoTildeSpecialForms}
Assume $E_Q$ is a regular $L_Q$-module with an admissible inner product.
For $P\le Q \in \Pl$, the ``restriction'' morphism $\widetilde k_{PQ}$
restricts to a natural morphism
\begin{equation*}
\widetilde k_{PQ}:\widetilde{\A}_{(2),\textup{sp}}(\Xhat_Q;\EE_Q)
\longrightarrow
\ihat_{P*}\Adec(\Xhat_P; \widetilde{\A}_{(2),\textup{sp}}(\bar{\mathscr
  A}_P^Q; \HH(\n_P^Q; E_Q), \wgt_P)_\infty)
\end{equation*}
such that $\i_P^* (\widetilde k_{PQ})$ is a quasi-isomorphism.
\end{enumerate}
\end{prop}

\begin{rem}
Note that in \itemref{itemRestrictionTildeSpecialForms} we can conclude
that $\ihat_P^*(\widetilde k_{PQ})$ is a quasi-isomorphism, but in
\itemref{itemRestrictionLTwoTildeSpecialForms} we only have that
$\i_P^*(\widetilde k_{PQ})$ is a quasi-isomorphism
\end{rem}

\begin{proof}
Part ~\itemref{itemTildeSpecialQuasiIsomorphism} follows from
\eqref{eqnSpecialFormsResolution} and Proposition
~\ref{propQuasiSpecialQuasiIsomorphism}%
\itemref{itemQuasiSpecialQuasiIsomorphism}.  In a special neighborhood $V=
\nil(\bar{\mathscr A}_P^Q(b) \times O_P)$ of a point $x\in X_P$ the
morphism $\widetilde k_{PQ}$ may be defined (in the notation of
\S\ref{ssectQuasiSpecialForms}) by
\begin{equation*}
\pr_P^*\psi\wedge (\pr^P)^*\o \longmapsto
\Adec(\Xhat_P;\lambda_P^Q)\circ k_{PQ}(\psi_\infty \otimes (\pr^P)^*\o),
\end{equation*}
where $\psi\in \Atilde\sp(\bar{\mathscr A}_P^Q(b);\EE_Q)$ and $\o \in
A(O_P)$ is $N_P^Q(\RR)$-invariant.  This indeed takes quasi-special forms
to quasi-special forms as may be verified using \eqref{eqnQuasiSpecialII}.
The assertion that $\ihat_P^* (\widetilde k_{PQ})$ is a quasi-isomorphism
follows from the corresponding assertion for $k_{PQ}$
\cite[Cor. 4.8]{refnSaperLModules} and the verification of
\eqref{eqnComposeRestrictionQuasiSpecialForms} may be left to the reader.
Part ~\itemref{itemRestrictionTildeSpecialForms} follows.  Furthermore the
$L^2$-norm near $x\in X_P$ on the left-hand side involves an integral over
the $\mathscr N_P^Q(\RR)'$-fibers of $\mathscr A_P^Q(b) \times O_P \to
\mathscr A_P^Q(b) \times \nil(O_P)$ which is not present in the right-hand
side.  Since the volume of the fiber over $(a\cdot x_0, x')$ is $\sim
\hsr_P(a)^{-2}$, where $a\in A_P^Q$, the integral over the fibers is
accounted for in the right-hand side by the weight $\wgt_P$.  Part
~\itemref{itemRestrictionLTwoTildeSpecialForms} follows.
\end{proof}

\section{Proof of Theorem ~\ref{thmLtwoLmodule}}
\label{sectProofThmLtwoLmodule}
We will prove Theorem ~\ref{thmLtwoLmodule} (from
\S\ref{sectLTwoCohomologyLModule}) by constructing morphisms of complexes
of sheaves
\begin{equation*}
\xymatrix{
 {\Ltwo(\Xhat;\EE)} & {\Sheaf_{\Xhat}(\Ltwo(E))} \ar[d]^t \\
{\widetilde{\A}_{(2),\textup{sp}}(\Xhat;\EE)} \ar[u]_r \ar[r]^s &
{\widetilde{\Sheaf}_{\Xhat}(\Ltwo(E))}
}
\end{equation*}
and show that they are all quasi-isomorphisms.

The complex $\Ltwo(\Xhat;\EE)$ is the usual $L^2$-cohomology sheaf as
defined in \S\ref{sectLTwoCohomology} and the sheaf
$\widetilde{\A}_{(2),\textup{sp}}(\Xhat;\EE)$ is the subcomplex obtained by
intersecting that with quasi-special forms as in
\S\ref{ssectQuasiSpecialForms}.  The morphism $r$ is the inclusion and it
is a quasi-isomorphism by Proposition
~\ref{propQuasiSpecialQuasiIsomorphism}%
\itemref{itemQuasiSpecialLtwoQuasiIsomorphism}.

The complex of sheaves $\Sheaf_{\Xhat}(\Ltwo(E))$ is the usual realization,
defined by applying \eqref{eqnRealizationFunctor} to
\eqref{eqnLTwoCohomologyLModule}, and the complex of sheaves
$\widetilde{\Sheaf}_{\Xhat}(\Ltwo(E))$ is defined similarly but using the
quasi-special differential forms functor.  Explicitly, the latter is
(compare \eqref{eqnRealizationFunctor})
\begin{equation}
  \left\{
  \begin{aligned}
    \widetilde{\Sheaf}_{\Xhat}(\Ltwo(E)) &= \bigoplus_{P\in \Pl}
    \ihat_{P*} \Adec(\Xhat_P;\EE_P),\\
    d_{\widetilde{\Sheaf}_{\Xhat}(\Ltwo(E))} &= \sum_{P\in\Pl} d_P +
    \sum_{P\le Q\in\Pl} \Adec(\Xhat_P;\widetilde f_{PQ})\circ \widetilde
    k_{PQ},
  \end{aligned}
  \right.
\label{eqnVariantRealizationFunctor}
\end{equation}
where $\widetilde k_{PQ}$ is from Proposition
~\ref{propRestrictionTildeSpecialForms}%
\itemref{itemRestrictionTildeSpecialForms} and $\widetilde f_{PQ}$ is
defined by
\begin{equation}
\label{eqnFTildeDefinition}
\widetilde f_{PQ}\colon \Adec(\bar{\mathscr A}_P^Q;
 \HH(\n_P^Q;E_Q))_\infty  \longrightarrow E_P[1], \qquad 
 v\otimes \psi
\longmapsto  f_{PQ}(v) \wedge (\pr^Q)^*\psi,
\end{equation}
for $v\in H(\n_P^Q;E_Q)$ and $\psi \in \Adec(\bar{\mathscr A}_P^Q)_\infty$.
The wedge product in \eqref{eqnFTildeDefinition} should be interpreted in
the following way.  By \eqref{eqnLambdaIsomorphism},
\eqref{eqnFRDefinition}, and \eqref{eqnLTwoCohomologyLModule}, the vector
$v$ is a sum of germs of forms on ${\mathscr A}_Q^G$, each of which satisfy
an $L^2$-condition with weight $\wgt_Q$ on ${\mathscr A}_R^G \times V$ for
various $R\ge Q$ and all $V\subseteq {\mathscr A}_Q^R$ relatively compact.
(For simplicity of notation we ignore here the finite dimensional
coefficients and the shift.)  By Lemma ~\ref{lemPullbackMorphism} the map
$f_{PQ}$ pulls these back to germs of forms on ${\mathscr A}_P^G$
satisfying a similar $L^2$-condition, except now $V\subseteq {\mathscr
A}_P^{R'}$ for some $R'\ge R$ and the weight is $\wgt_P$.  The wedge
product with $(\pr^Q)^*\psi$, the germ of a pullback of a smooth form on
${\mathscr A}_P^Q$, preserves this $L^2$-condition.

The morphism $t$ is the sum of the  inclusions
\begin{equation*}
\ihat_{P*}\Asp(\Xhat_P;\EE_P) \hookrightarrow \ihat_{P*}
\Adec(\Xhat_P;\EE_P).
\end{equation*}
The complexes $\Sheaf_{\Xhat}(\Ltwo(E))$ and
$\widetilde{\Sheaf}_{\Xhat}(\Ltwo(E))$ are filtered by the parabolic rank
of $P$; the morphism induced by $t$ on the associated graded complexes is a
quasi-isomorphism by Proposition ~\ref{propRestrictionTildeSpecialForms}%
\itemref{itemTildeSpecialQuasiIsomorphism} and hence so is $t$.

It remains to define $s$ and prove that it is a quasi-isomorphism.  For
each $P \in \Pl$, the morphism
\begin{equation*}
\widetilde k_{PG} \colon \widetilde{\A}_{(2),\textup{sp}}(\Xhat;\EE)
\longrightarrow \ihat_{P*} \Adec(\Xhat_P;\widetilde
\A_{(2),\textup{sp}}(\bar{\mathscr A}_P^G; \HH(\n_P;E), \wgt_P)_\infty)
\end{equation*}
from Proposition ~\ref{propRestrictionTildeSpecialForms}%
\itemref{itemRestrictionLTwoTildeSpecialForms} may be composed with the
morphism induced by
\begin{equation*}
\widetilde \A_{(2),\textup{sp}}(\bar{\mathscr A}_P^G; \HH(\n_P;E),
\wgt_P)_\infty
\subseteq \A_{(2)}(\bar{\mathscr A}_P^G; \HH(\n_P;E), \wgt_P)_\infty = 
F_{\{P\}} \subseteq E_{P}.
\end{equation*}
We can thus define
\begin{equation}
\label{eqnSMorphismDefinition}
s = \sum_P \widetilde k_{PG} \colon
\widetilde{\A}_{(2),\textup{sp}}(\Xhat;\EE) \longrightarrow
\bigoplus_{P\in \Pl} \ihat_{P*} \Adec(\Xhat_P;\EE_P) = 
\widetilde{\Sheaf}_{\Xhat}(\Ltwo(E)).
\end{equation}

{}To prove that $s$ is a quasi-isomorphism we will fix $P\in \Pl$ and prove
that $\i_P^*s$ is a quasi-isomorphism.  Consider the morphisms
\begin{equation}
\label{eqnAnalyzeS}
\raisebox{.3\totalheight}{\makebox[0pt][c]{\xymatrix{
*+!{\i_P^*\widetilde{\A}_{(2),\textup{sp}}(\Xhat;\EE)}
\ar[r]^-{\sum_{R\ge P} \i_P^*(\widetilde k_{RG})}
\ar[d]^-{\i_P^* (\widetilde k_{PG})} &
*+!{\smash{\bigoplus_{R\ge P}}\i_P^*
  \Adec(\Xhat_R;\EE_R)} \ar[d]^-{\sum_{R\ge P} \i_P^*(\widetilde k_{PR})}\\
*+!{\A(X_P; \widetilde{\A}_{(2),\textup{sp}}(\bar{\mathscr A}_P^G;
  \HH(\n_P;E), \wgt_P)_\infty)} &
*+!{\bigoplus_{R\ge P} \A(X_P; \Adec(\bar{\mathscr
    A}_P^R;\HH(\n_P^R;E_R))_\infty)\  ,}
}}}
\end{equation}
where the upper right-hand group has the differential
\begin{equation}
\label{eqnUpperDifferential}
\sum_{R\ge P} d_R +
\sum_{S\ge R\ge P} \Adec(\Xhat_R;\widetilde f_{RS})\circ \widetilde k_{RS}
\end{equation}
and the bottom right-hand group has the differential
\begin{equation}
\label{eqnLowerDifferential}
\sum_{R\ge P} (d_P + d_{\mathscr A_P^R}) +
\sum_{S\ge R\ge P} \A(X_P; \Adec(\bar{\mathscr A}_P^R;
H(\n_P^R;\widetilde f_{RS}))_\infty) \circ 
\A(X_P;\widetilde m_{RS}).
\end{equation}
The top horizontal map then represents $\i_P^* s$.  Both the left and right
vertical morphisms are quasi-isomorphisms by Proposition
~\ref{propRestrictionTildeSpecialForms}, parts
\itemref{itemRestrictionLTwoTildeSpecialForms} and
\itemref{itemRestrictionTildeSpecialForms} respectively.  (On the
right-hand side, one applies the proposition to the graded morphism
associated to the double filtration by parabolic rank of $R$ and by degree
of $E_R$.)

We can complete \eqref{eqnAnalyzeS} to a commutative diagram by defining
the morphism 
\begin{equation}
\label{eqnBottomHorizontal}
\A(X_P; \widetilde{\A}_{(2),\textup{sp}}(\bar{\mathscr A}_P^G;
  \HH(\n_P;E), \wgt_P)_\infty) \longrightarrow
\bigoplus_{R\ge P} \A(X_P; \Adec(\bar{\mathscr
    A}_P^R;\HH(\n_P^R;E_R))_\infty)
\end{equation}
to be $\sum_{R\ge P} \A(X_P;\widetilde m_{RG})$.  More precisely, each term
in this sum represents the morphism induced by
\begin{equation*}
\widetilde{\A}_{(2),\textup{sp}}(\bar{\mathscr A}_P^G;
  \HH(\n_P;E), \wgt_P)_\infty \longrightarrow
\Adec(\bar{\mathscr A}_P^R;
  \widetilde{\A}_{(2),\textup{sp}}(\bar{\mathscr A}_R^G;
  \HH(\n_P;E), \wgt_R)_\infty)_\infty
\end{equation*}
from \eqref{eqnRestrictionSplitLTwoQuasiSpecialForms}, followed by the
morphisms induced by
\begin{equation*}
\widetilde{\A}_{(2),\textup{sp}}(\bar{\mathscr A}_R^G;
  \HH(\n_P;E), \wgt_R)_\infty \longtildearrow
H(\n_P^R; \widetilde{\A}_{(2),\textup{sp}}(\bar{\mathscr A}_R^G;
  \HH(\n_R;E), \wgt_R)_\infty)
\end{equation*}
from \eqref{eqnLambdaIsomorphism} and by
\begin{equation*}
 \widetilde{\A}_{(2),\textup{sp}}(\bar{\mathscr A}_R^G;
  \HH(\n_R;E), \wgt_R)_\infty \subseteq \A_{(2)}(\bar{\mathscr
  A}_R^G; \HH(\n_R;E), \wgt_R)_\infty =
  F_{\{R\}} \subseteq E_R.
\end{equation*}

We need to show that \eqref{eqnBottomHorizontal} is a quasi-isomorphism.
The proof is parallel to that of Proposition ~\ref{propLocalLtwo} of
\S\ref{sectLocalLtwo} so we will be brief.  Define filtrations on the two
complexes: on the left-hand side of \eqref{eqnBottomHorizontal} use the
trivial filtration in which $F^0$ is the entire complex and $F^p=0$ for
$p>0$; on the right-hand side use \eqref{eqnLTwoCohomologyLModule} to
re-express the sum over $R\ge P$ as a sum over $\R$ with $R_0\ge P$ and let
$F^p$ consist of terms such that $\#(\R\setminus\{P\})\ge p$.

For $p=0$ the graded morphism associated to \eqref{eqnBottomHorizontal} is
the inclusion
\begin{equation*}
\A(X_P; \widetilde{\A}_{(2),\textup{sp}}(\bar{\mathscr A}_P^G;
  \HH(\n_P;E), \wgt_P)_\infty) \longrightarrow
\A(X_P;\A_{(2)}(\bar{\mathscr A}_{P}^G;
\HH(\n_{P};E), \wgt_P)_\infty)
\end{equation*}
which is a quasi-isomorphism by Proposition
~\ref{propSplitQuasiSpecialQuasiIsomorphism}%
\itemref{itemSplitQuasiSpecialLtwoQuasiIsomorphism}.  For $p>0$ it is the
map of $0$ into
\begin{equation*}
\bigoplus_{\substack{\R \mid R_0 > P \\ \#\R = p}}
  \left( \A(X_P;\Adec(\bar{\mathscr A}_P^{R_0};\HH(\n_P^{R_0};F_{\R}))_\infty)
  \oplus \A(X_P;F_{\{P\}\cup\R}) \right).
\end{equation*}
This is a direct sum of complexes so it suffices to show that the summand
for a given $\R$ is acyclic.  This summand is a shifted mapping cone for
$\A(X_P;\widetilde f_{PR_0})$:
\begin{equation*}
\xymatrix @C=-.1in{
{\A(X_P;\Adec(\bar{\mathscr A}_P^{R_0};\HH(\n_P^{R_0};F_{\R}))_\infty)}
\ar[rr]^-{\A(X_P;\widetilde f_{PR_0})}  & {} &
{\A(X_P;F_{\{P\}\cup\R})[1]} \\
{} & {\A(X_P;\HH(\n_P^{R_0};F_{\R}))_\infty\rlap{\ .}} \ar[lu]
\ar[ru]_*!<7pt,0pt>{\labelstyle-\A(X_P;u_{\{P\}\cup\R,\R})}
}
\end{equation*}
Since Lemma ~\ref{lemPullbackMorphismQuasiIsomorphism} implies that
$u_{\{P\}\cup\R,\R}$ is a quasi-isomorphism and Proposition
~\ref{propSplitQuasiSpecialQuasiIsomorphism}%
\itemref{itemSplitQuasiSpecialQuasiIsomorphism} together with the usual
Poincar\'e lemma implies that the left diagonal map is a quasi-isomorphism,
the proof is complete. \qed

\section{The Micro-support of the $L^2$-cohomology $\L$-module}
\label{sectMicroSupportLtwo}
Let $E$ be a regular $G$-module and let $\Ltwo(E)$ be the corresponding
$\L$-module as in \S\ref{sectLTwoCohomologyLModule}.  The calculation of
the micro-support $\mS(\Ltwo(E))$ that follows is similar to that for
weighted cohomology in \cite[\S16]{refnSaperLModules}, although more
complicated by the possible presence of infinite-dimensional local
cohomology; I am grateful to an anonymous referee of
\cite{refnSaperLModules} for comments that simplified the proof.

\subsection{}
It is helpful to use the language of lattices.  Let $P\in \Pl$.  The
partially ordered set $[P,G]= \{\, Q\in\Pl \mid P\le Q\le G\,\}$ is a
Boolean lattice; in fact the map $Q\mapsto \D_P^Q$ is an isomorphism of
$[P,G]$ onto the lattice of subsets of $\D_P$.  Given $Q\in [P,G]$, let
$(P,Q)\in [P,G]$ denote the complementary element determined by
$\D_P^{(P,Q)} = \D_P\setminus \D_P^Q$.  For $Q$, $R\in [P,G]$, let $Q\vee
R$ and $Q\wedge R$ denote the elements of $[P,G]$ determined by $\D_P^Q\cup
\D_P^R$ and $\D_P^Q\cap \D_P^R$ respectively.  If $R\le S$, the subset
$[R,S] = \{\, Q\in\Pl \mid R\le Q\le S\,\}$ is a Boolean sublattice.

\subsection{}
Let $P\in \Pl$ and let $V$ be an irreducible $L_P$-module.  Recall Zucker's
calculation of the $L^2$-cohomology of $\mathscr A_P^G(b)$
\cite[(4.51)]{refnZuckerWarped}:
\begin{lem}
\label{lemSplitLtwoCohomologyCalculation}
For any $b\in \mathscr A_P^G$,
\begin{equation*}
H_{(2)}(\mathscr A_{P}^G(b);\VV,\wgt_P) =
\begin{cases}
\CC & \text{if $\langle \xi_V + \hsr_P, \b_\al\spcheck\rangle > 0$ for all $\al
  \in \D_P$,} \\
0 & \text{if $\langle \xi_V + \hsr_P, \b_\al\spcheck\rangle < 0$ for any  $\al
  \in \D_P$.}
\end{cases}
\end{equation*}
In the remaining case where $\langle \xi_V + \hsr_P, \b_\al\spcheck\rangle
\ge 0$ for all $\al \in \D_P$ and $r = \#\{\,\al \mid \langle \xi_V +
\hsr_P, \b_\al\spcheck \rangle = 0\,\}>0$, then $H^i_{(2)}(\mathscr
A_{P}^G(b);\VV,\wgt_P)$ is nonzero only if $i\in [1, r]$ in which case it
is infinite dimensional.
\end{lem}

Define $S_V$,
$T_V$, $T_V'\ge P$ by
\begin{align*}
\D_P^{S_V} &= \{\, \al\in\D_P\mid \langle \xi_V + \hsr_P, \b_\al\spcheck
\rangle > 0\,\}, \\
\D_P^{T_V} &= \{\, \al\in\D_P\mid \langle \xi_V + \hsr_P, \b_\al\spcheck
\rangle < 0\,\}, \\
\D_P^{T_V'}& = \{\, \al\in\D_P\mid \langle \xi_V + \hsr_P, \b_\al\spcheck
\rangle \le 0\,\}.
\end{align*}

\begin{lem}
\label{lemLocalLtwoCohomology}
\begin{equation*}
H(\i_P^* \Ltwo(E))_V \cong
\begin{cases}
H_{(2)}(\mathscr A_{S_V}^G(b);\CC)\otimes H(\n_P;E)_V & \text{if $T_V= P$,}
\\
0 & \text{otherwise,}
\end{cases}
\end{equation*}
for any $b \in \mathscr A_{S_V}^G$.  If $T_V=P$, the group
$H^i_{(2)}(\mathscr A_{S_V}^G(b);\CC)$ is $\CC$ if $S_V=G$ and otherwise is
nonzero \textup(and infinite dimensional\textup) only for degrees in $[1,
\dim A_{S_V}^G]$.
\end{lem}
\begin{proof}
Proposition ~\ref{propLocalLtwo} implies that $ H(\i_P^* \Ltwo(E))_V \cong
H_{(2)}(\mathscr A_P^G(b);\VV,\wgt_P)\otimes H(\n_P;E)_V$ which by Lemma
~\ref{lemSplitLtwoCohomologyCalculation} is zero unless $T_V=P$.  Zucker's
K\"unneth formula \cite[(2.34)(i)]{refnZuckerWarped} applied to the
decomposition $\mathscr A_P^G(b) \cong \mathscr A_{S_V}^G(b) \times
\mathscr A_{T_V'}^G(1)$ from \cite[4.3(3)]{refnBorelSerre} (which is
different from \eqref{eqnSplitDecomposition}) yields
\begin{equation*}
H_{(2)}(\mathscr A_P^G(b); \VV,\wgt_P) \cong H_{(2)}(\mathscr
A_{S_V}^G(b);\VV, \wgt_{S_V}) \otimes  H_{(2)}(\mathscr
A_{T_V'}^G(1);\VV, \wgt_{T_V'})
\end{equation*}
provided one of the factors is finite dimensional.  Lemma
~\ref{lemSplitLtwoCohomologyCalculation} implies the second factor is $\CC$
and the remaining assertions.
\end{proof}

\begin{lem}
\label{lemHNilpotentLocalLtwoCohomology}
For $P\le R$,
\begin{equation*}
H(\n_P^R;H(\i_R^* \Ltwo(E)))_V \cong
\begin{cases}
H_{(2)}(\mathscr A_{R\vee S_V}^G(b);\CC)\otimes H(\n_P;E)_V & \text{if $T_V
  \le R$,} \\
0 & \text{otherwise.}
\end{cases}
\end{equation*}
\end{lem}
\begin{proof}
If the group is nonzero then $H(\n_P^R;W)_V \neq 0$ for some irreducible
regular $L_R$-module $W$ occuring in $H(\i_R^* \Ltwo(E))$.  Let $\lambda$
be the highest weight of $W$ with respect to a Cartan subalgebra of
$\levi_R$ and a positive system of roots which contains those in $\n_P^R$.
Under these conditions, Kostant's theorem \cite[Theorem ~5.14]{refnKostant}
implies that $H(\n_P^R;W)_V =V$ and the highest weight of $V$ is $w(\lambda
+ \hsr)- \hsr$, where $\hsr$ is one-half the sum of the positive roots of
$\levi_R$ and $w$ belongs to a certain subset of the Weyl group of $\levi_R$.
Thus $W$ is uniquely determined and $\xi_V|_{\sa_R} = \xi_W$.  The lemma
now follows from Lemma ~\ref{lemLocalLtwoCohomology} (applied to
$H(\i_R^*\Ltwo(E))_W$), since it is clear that $T_W = R\vee T_V$ and $S_W =
R\vee S_V$.
\end{proof}

\begin{lem}
\label{lemIsomorphismLocalLtwoCohomology}
Suppose that $P\le R\le R'$ such that $R\vee S_V = R'\vee S_V$.  Then the
natural morphism
\begin{equation*}
H(\n_P^R;H(\i_R^* \Ltwo(E)))_V \longrightarrow  H(\n_P^{R'};H(\i_{R'}^*
\Ltwo(E)))_V
\end{equation*}
corresponds \textup(up to sign\textup) to the identity under the
isomorphisms of Lemma ~\textup{\ref{lemHNilpotentLocalLtwoCohomology}}.
\end{lem}
\begin{proof}
The natural morphism is given by \eqref{eqnRetrictionLtwo} in Corollary
~\ref{corNaturalMorphismLtwo} (aside from the application of
$H(\n_P^R;\cdot)$).  The lemma follows easily.
\end{proof}

\subsection{}
\label{ssectLatticeSplitCohomology}
Consider $S_1\le S_2$ and order $\D_{S_1} = \{\al_1,\dots,\al_r\}$.  Denote
by $A_{(2),R}(\mathscr A_{S_1}^G(b);\CC)$ the sections of
$\A_{(2),R}(\bar{\mathscr A}_{S_1}^G;\CC)$ over $\bar{\mathscr
A}_{S_1}^G(b)$.  Define the double complex
\begin{equation*}
A_{(2)}(\mathscr A_{[S_1,S_2]}^G(b);\CC) = \bigoplus_{S_1 \le R\le
S_2} A_{(2),R}(\mathscr A_{S_1}^G(b);\CC)[-\#\D_{S_1}^R]
\end{equation*}
where the horizontal differential between the $R'$ and $R$ terms (when
$\D_{S_1}^R = \D_{S_1}^{R'}\cup\{\al_i\}$) is $(-1)^i g_{S_1R,S_1R'}$.
Since $H(A_{(2),R}(\mathscr A_{S_1}^G(b);\CC)) \cong H_{(2)}(\mathscr
A_R^G(b);\CC)$ by Lemma ~\ref{lemPullbackMorphismQuasiIsomorphism}, the
$E_1$-term of the spectral sequence for the total complex is
\begin{equation}
\label{eqnIntervalCohomologySpectralSequence}
\bigoplus_{S_1\le R\le S_2} H_{(2)}(\mathscr A_R^G(b);\CC)[-\#\D_{S_1}^R],
\end{equation}
where the terms of $d_1$ are given by Corollary
~\ref{corNaturalMorphismLtwo}.  Note that all of these terms are infinite
dimensional with the exception of the $R=S_2$ term in the case that
$S_2=G$.  We denote the cohomology of the total complex by
$H_{(2)}(\mathscr A_{[S_1,S_2]}^G(b);\CC)$.  We may similarly define a
complex and cohomology for open and half-open intervals such as
$(S_1,S_2]$.

The cohomology $H_{(2)}(\mathscr A_{[S_1,S_2]}^G(b);\CC)$ is always
nonzero; for example, the spectral sequence
\eqref{eqnIntervalCohomologySpectralSequence} shows that it does not vanish
in degree $\dim \mathscr A_{S_1}^G$.  Furthermore, unless $S_1=S_2=G$, the
cohomology is infinite dimensional.

\begin{prop}
\label{propLTwoCohomologyWithSupports}
For $P\le Q\in \Pl$ and an irreducible $L_P$-module $V$,
\begin{multline*}
H(\i_P^*\ihat_Q^!\Ltwo(E))_V \cong \\
\begin{cases}
  H_{(2)}(\mathscr A_{[T_V\vee S_V,(P,Q)\vee S_V]}^G(b);\CC)
  \otimes H(\n_P;E)_V[-\#\D_P^{T_V}] & \text{if $(P,T_V') \le Q \le
  (P,T_V)$,} \\
  0 & \text{otherwise.}
				   \end{cases}
\end{multline*}
\end{prop}

\begin{proof}
Consider the short exact sequence
\begin{equation}
\label{eqnLTwoShortExactSequence}
0 \longrightarrow \i_P^*\ihat_Q^!\Ltwo(E)_V \longrightarrow \i_P^*\Ltwo(E)_V
\longrightarrow \i_P^* \jhat_{Q*}\jhat_Q^* \Ltwo(E)_V \longrightarrow 0
\end{equation}
obtained from \eqref{eqnShortExactSequenceWithSupport} by taking the
$V$-isotypical component.  There is a Mayer-Vietoris spectral sequence
\cite[Lemma ~3.7]{refnSaperLModules} abutting to $H( \i_P^*
\jhat_{Q*}\jhat_Q^* \Ltwo(E))_V$ with
\begin{equation}
\label{eqnLinkRelativeMayerVietoris}
E_1^{p,\cdot} = \bigoplus_{\substack{P < R \le (P,Q) \\
    \#\D_P^{R} = p+1}} H(\n_P^{R}; H(\i_{R}^* \Ltwo(E)))_V.
\end{equation}
By Lemma ~\ref{lemHNilpotentLocalLtwoCohomology} the term indexed by $R$
will vanish unless $T_V\le R$.  Thus \eqref{eqnLinkRelativeMayerVietoris}
vanishes unless $T_V\le (P,Q)$, and in this case the terms of
\eqref{eqnLinkRelativeMayerVietoris} are indexed by $R\in [T_V, (P,Q)]$
(with $P$ excluded if $T_V=P$).

In the lattice $[T_V,G]$, the complement of $T_V'$ is $T_V\vee S_V$.  Any
element of $[T_V,G]$ may thus be expressed as the join of its intersection
with $T_V'$ and its intersection with $T_V\vee S_V$.  We apply this to
elements of the sublattice $[T_V, (P,Q)]$; if we fix the first intersection
to be $\tilde T_V$ and let the second intersection vary, we obtain elements
$R$ on the dotted line below:
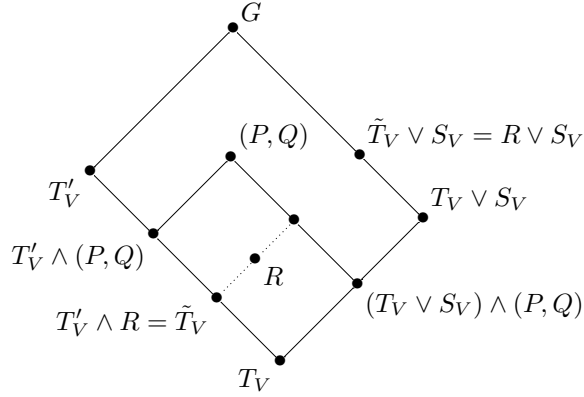
\begin{equation*}
\vcenter{%
\entrymodifiers={[o]}%
\xymatrix @ur @R=2.3pc @C=1.2pc {
{\bullet}\save[]*+!RU{T_V'}\restore \ar@{-}[rrrr] \ar@{-}[d] &
{} &
{} &
{} &
{\bullet}\save[]*+!LD{G}\restore \ar@{-}[dd] \\
{\bullet}\save[]*+!RU{T_V'\wedge (P,Q)}\restore \ar@{-}[rr] \ar@{-}[d] &
{} &
{\bullet}\save[]*+!LD{(P,Q)}\restore \ar@{-}[d] &
{} &
{} \\
{\bullet}\save[]*+!RU{T_V'\wedge R = \tilde T_V}\restore \ar@{.}[r]
\ar@{-}[d] &
{\bullet}\save[]*+!LU{R}\restore \ar@{.}[r] &
{\bullet}\save[]*+!LD{}\restore \ar@{-}[d] &
{} &
{\bullet}\save[]*+!LD{\tilde T_V\vee S_V = R\vee S_V}\restore \ar@{-}[d] \\
{\bullet}\save[]*+!RU{T_V}\restore \ar@{-}[rr] &
{} &
{\bullet}\save[]*+!UL{(T_V\vee S_V)\wedge(P,Q)}\restore \ar@{-}[rr] &
{} &
{\bullet}\save[]*+!LD{T_V\vee S_V}\restore
}}
\end{equation*}
If we then vary $\tilde T_V$ we obtain a decomposition
\begin{equation}
\label{eqnLatticeDecomposition}
[T_V, (P,Q)] = \coprod_{\tilde T_V \in [T_V, T_V'\wedge (P,Q)]}  [\tilde
  T_V, (\tilde T_V \vee S_V) \wedge (P,Q)]
\end{equation}
The elements $R$ in a fixed component of \eqref{eqnLatticeDecomposition}
(say indexed by $\tilde T_V$) will all have the same value of $R\vee S_V$,
namely $\tilde T_V \vee S_V$, and the same value of $T_V'\wedge R$, namely
$\tilde T_V$.  The corresponding terms of
\eqref{eqnLinkRelativeMayerVietoris} will then all be isomorphic by
Lemma ~\ref{lemHNilpotentLocalLtwoCohomology}.

In view of the preceding discussion, filter the complex $(E_1,d_1)$ by
$\#\D_P^{T_V'\wedge R}$.  The associated graded complex is a direct sum of
complexes indexed by $\tilde T_V \in [T_V, T_V'\wedge (P,Q)]$.  The complex
associated to a given $\tilde T_V$ is
\begin{equation}
\label{eqnTildeTComplex}
\bigoplus_{\substack{R\in [\tilde T_V, (\tilde T_V \vee S_V) \wedge
(P,Q)]\\ R\neq P}}
H_{(2)}(\mathscr A_{\tilde T_V\vee S_V}^G(b);\CC)\otimes H(\n_P;E)_V
[-\#\D_P^R +1]
\end{equation}
with differential $\sum_{\#\D_P^{R'} = \#\D_P^R+1} \pm\id_{R,R'}$ by Lemma
~\ref{lemIsomorphismLocalLtwoCohomology}; here $\id_{R,R'}$ denotes the
identity morphism between the $R$-term and the $R'$-term.

First assume $P < T_V$.  Then Lemma ~\ref{lemLocalLtwoCohomology} and the
long exact sequence associated to \eqref{eqnLTwoShortExactSequence} imply
that $H(\i_P^*\ihat_Q^!\Ltwo(E))_V \cong H(\i_P^* \jhat_{Q*}\jhat_Q^*
\Ltwo(E))_V[-1]$.  If furthermore $(P,Q)\notin [T_V, T_V']$, then the
cohomology of \eqref{eqnTildeTComplex} vanishes: aside from a shift it is
the simplicial cohomology of the cone over the simplex with vertices $\D_{
\tilde T_V}^{(\tilde T_V \vee S_V) \wedge (P,Q)}$.  On the other hand, if
$(P,Q)\in [T_V, T_V']$, that is, $(P,T_V') \le Q \le (P,T_V)$, then
\eqref{eqnTildeTComplex} reduces to
\begin{equation*}
H_{(2)}(\mathscr A_{\tilde T_V\vee S_V}^G(b);\CC)\otimes H(\n_P;E)_V
[-\#\D_P^{\tilde T_V} +1]
\end{equation*}
and the spectral sequence \eqref{eqnLinkRelativeMayerVietoris} is
isomorphic to the spectral sequence
\eqref{eqnIntervalCohomologySpectralSequence} for $H_{(2)}(\mathscr
A_{[T_V\vee S_V,(P,Q)\vee S_V]}^G(b);\CC)$ (tensored with $H(\n_P;E)_V$ and
shifted by $1-\#\D_P^{T_V}$).

On the other hand, assume $P=T_V$.  Now if $(P,Q)\notin [T_V, T_V']$, the
cohomology of \eqref{eqnTildeTComplex} vanishes except for the case $\tilde
T_V= P$, in which case it is $H_{(2)}(\mathscr A_{S_V}^G(b);\CC)\otimes
H(\n_P;E)_V$.  Thus the spectral sequence
\eqref{eqnLinkRelativeMayerVietoris} degenerates and its cohomology will be
canceled in the long exact sequence associated to
\eqref{eqnLTwoShortExactSequence} by $H(\i_P^*\Ltwo(E))_V$ (use Lemma
~\ref{lemLocalLtwoCohomology}).  If $(P,Q)\in [T_V, T_V']$, the above
argument shows that \eqref{eqnLinkRelativeMayerVietoris} abuts to
$H_{(2)}(\mathscr A_{(S_V,(P,Q)\vee S_V]}^G(b);\CC) \otimes
H(\n_P;E)_V[1]$.  We have a commutative diagram with exact rows
\begin{equation*}
\xymatrix @1 @C-11pt{
{} \ar[r]  &  {H^{i-1}(\i_P^* \jhat_{Q*}\jhat_Q^* \Ltwo(E))_V}
\ar[r] \ar[d]  & {H^i(\i_P^*\ihat_Q^!\Ltwo(E))_V} \ar[r] \ar[d]
 & {H^i(\i_P^*\Ltwo(E))_V} \ar[r] \ar[d] & {} \\
{} \ar[r] & {H^i_{(2)}(\mathscr A_{(S_V,(P,Q)\vee S_V]}^G(b);\CC)} \ar[r] &
{H^i_{(2)}(\mathscr A_{[S_V,(P,Q)\vee S_V]}^G(b);\CC)} \ar[r]  &
{H^i_{(2)}(\mathscr A_{S_V}^G(b);\CC)} \ar[r] & {} 
}
\end{equation*}
where to save space we have omitted the tensor product with $H(\n_P;E)_V$
in the bottom row.  We have already noted that the first vertical arrow is
an isomorphism, while the last vertical arrow is an isomorphism by Lemma
~\ref{lemLocalLtwoCohomology}.  The proposition now follows from the
$5$-lemma.
\end{proof}

\begin{thm}
\label{thmLTwoMicroSupport}
For a regular $G$-module $E$, the micro-support $\mS(\Ltwo(E))$ consists of
those irreducible $L_P$-modules $V$ satisfying
\begin{enumerate}
\item\label{itemNilpotentType} $H(\n_P;E)_V\neq 0$,
\item\label{itemConjugateSelfContragredient} $(V|_{\lsp0L_P})^* \cong
\overline{V|_{\lsp0L_P}}$, and
\item\label{itemSplitVanishing} $(\xi_V + \hsr_P)|_{\sa_P^G} =0$.
\end{enumerate}
For such a $V$ and any $Q\in [P,G]$,
\begin{equation*}
\Type_{Q,V}(\Ltwo(E)) = H_{(2)}(\mathscr A_{[P,(P,Q)]}^G(b);\CC) \otimes
H(\n_P;E)_V.
\end{equation*}
The weak micro-support $\mS_w(\Ltwo(E))$ is similarly characterized by
omitting condition \itemref{itemConjugateSelfContragredient}.
\end{thm}

\begin{proof}
Let $V$ be an irreducible $L_P$-module.  By the definition of micro-support
in \S\ref{sectMicroSupport} and Proposition
~\ref{propLTwoCohomologyWithSupports}, $V\in \mS(\Ltwo(E))$ if and only if
conditions \itemref{itemNilpotentType},
\itemref{itemConjugateSelfContragredient}, and
\begin{equation}
\label{eqnRangeIntersection}
[Q_V, Q'_V] \cap [(P,T_V') , (P,T_V)] \neq \emptyset
\tag*{\itemref{itemSplitVanishing}$'$}
\end{equation}
are satisfied.  Clearly \itemref{itemSplitVanishing} implies
\ref{eqnRangeIntersection}.  Conversely we will assume
\ref{eqnRangeIntersection} holds and prove \itemref{itemSplitVanishing};
together with Proposition ~ \ref{propLTwoCohomologyWithSupports} this will
establish the theorem.

By Langlands's ``geometric lemma'' \cite[IV,
\S6.11]{refnBorelWallachSecondEdition}, there exists one and only one $R\in
[P,G]$ such that
\begin{subequations}
\begin{align}
\label{eqnPartitionPieceA}
&\langle \xi_V + \hsr_P,
\b_\al^R{}\spcheck \rangle \ge 0 &&\text{for $\al \in \D_P^R$, and}\\
\label{eqnPartitionPieceB}
&\langle \xi_V + \hsr_P, \g\spcheck_R
\rangle < 0 &&\text{for $\g\in \D_P \setminus \D_P^R$.}
\end{align}
\end{subequations}
Since $\g\spcheck{}^R \in -\cl{\sa_P^{R*+}} $ for $\g \in \D_P \setminus
\D_P^R$, \eqref{eqnPartitionPieceA} implies that $\langle \xi_V + \hsr_P,
\g\spcheck{}^R \rangle \le 0$.  Together with \eqref{eqnPartitionPieceB},
this yields $\langle \xi_V + \hsr_P, \g\spcheck{} \rangle < 0$ for all $\g
\in \D_P \setminus \D_P^R$, that is, $ (P,R) \le Q_V$.  However
\ref{eqnRangeIntersection} implies that $Q_V \le (P,T_V)$.  We conclude
that $\langle \xi_V + \hsr_P, \b_\g\spcheck \rangle \ge 0$ for all $\g\in
\D_P \setminus \D_P^R$, which means that $(\xi_V + \hsr_P)_R \in
\cl{\lsp+\sa_R^*}$.  However \eqref{eqnPartitionPieceB} implies that
$(\xi_V + \hsr_P)_R \in - \sa_R^{*+}$.  Since $\cl{\lsp+\sa_R^*} \cap (
-\sa_R^{*+}) = \emptyset$ unless $R=G$, we see that
\begin{equation}
\label{eqnRootCone}
\xi_V + \hsr_P \in \cl{\lsp+\sa_P^*}.
\end{equation}

Similarly there exists a unique $R\in [P,G]$ such that
\begin{subequations}
\begin{align}
\label{eqnPrimePartitionPieceA}
&\langle \xi_V + \hsr_P,
\b_\al^R{}\spcheck \rangle > 0 &&\text{for $\al \in \D_P^R$, and} \\
\label{eqnPrimePartitionPieceB}
&\langle \xi_V + \hsr_P, \g\spcheck_R
\rangle \le 0 &&\text{for $\g\in \D_P \setminus \D_P^R$.}
\end{align}
\end{subequations}
Since $\b_{\al R}\spcheck \in \cl{\sa_R^{*+}}$ for $\al \in \D_P^R$,
equation \eqref{eqnRootCone} implies that $\langle \xi_V + \hsr_P, \b_{\al
R}\spcheck \rangle \ge 0$.  Together with \eqref{eqnPrimePartitionPieceA},
this yields $\langle \xi_V + \hsr_P, \b_\al\spcheck \rangle > 0$ for all
$\al \in \D_P^R$, that is, $R \le (P,T_V')$.  However
\ref{eqnRangeIntersection} implies that $(P,T_V') \le Q_V'$.  We conclude
that $\langle \xi_V + \hsr_P, \al\spcheck \rangle \le 0$ for all
$\al\in\D_P^R$, which means that $(\xi_V + \hsr_P)^R \in
-\cl{\sa_P^{R*+}}$.  However \eqref{eqnPrimePartitionPieceA} implies that
$(\xi_V + \hsr_P)^R \in \lsp+\sa_P^{R*}$.  Since $(-\cl{\sa_P^{R*+}}) \cap
\lsp+\sa_P^{R*} = \emptyset$ unless $R=P$, we see that $\xi_V + \hsr_P \in
-\cl{\sa_P^{*+}}$.  This together with \eqref{eqnRootCone} establishes
\itemref{itemSplitVanishing} since $ \cl{\lsp+\sa_P^*} \cap (-
\cl{\sa_P^{*+}})= \sa_G^*$.
\end{proof}

A parabolic $\RR$-subgroup $P_0$ of $G$ is called \emph{fundamental} if
$\p_0$ contains a fundamental (that is, maximally compact) Cartan
subalgebra of $\mathfrak g$.
\begin{cor}
\label{corMicroSupportFundamentalCharacterization}
For $P\in \Pl$, there exists an irreducible $L_P$-module $V \in
\mS(\Ltwo(E))$ if and only if $(E|_{\lsp0G})^* \cong
\overline{E|_{\lsp0G}}$ and $P$ contains a fundamental parabolic
$\RR$-subgroup of $G$.  The type of $V$ is finite dimensional if and only
if $P=G$.
\end{cor}

\begin{proof}
Apply Theorem ~\ref{thmLTwoMicroSupport},
\cite[3.6(iii)(iv)]{refnBorelCasselman} (see also \cite[Lemma
~8.8]{refnSaperLModules}) and the remark at the end of
\S\ref{ssectLatticeSplitCohomology}.
\end{proof}

Recall that a symmetric space $D=G(\RR)/KA_G$ is called \emph{equal-rank}
if $\CCrank G = \rank K + \rank A_G$.  Any Hermitian symmetric space is
equal-rank and every equal-rank symmetric space has even dimension.  The
symmetric spaces associated to $G(\RR)=\SO(2p,2q+1)$ where $p>1$ are
examples of non-Hermitian equal-rank symmetric spaces.

\begin{cor}
\label{corLtwoMicropurity}
If $X$ is an arithmetic quotient of an equal-rank symmetric space and $E$
is an irreducible regular $G$-module, then $\mS(\Ltwo(E)) = \{E\}$.
\end{cor}

\begin{proof}
Since $D$ is equal-rank, the only fundamental parabolic $\RR$-subgroup of
$G$ is $G$ itself.  Furthermore $(E|_{\lsp0G})^* \cong
\overline{E|_{\lsp0G}}$ for any $G$-module $E$ \cite[1.5,
1.6]{refnBorelCasselman}.  Now apply Corollary
~\ref{corMicroSupportFundamentalCharacterization}.
\end{proof}

\section{The Conjectures of Borel and Zucker}
\label{sectBorelZuckerConjecture}
Associated to any finite-dimensional irreducible representation $\sigma$ of
$G(\RR)$, Satake \cite{refnSatakeCompact} constructs a compactification
$\Dstar_\sigma$ of $D$ which is a disjoint union of so-called \emph{real
boundary components}; $D$ is always a real boundary component and the
others are symmetric spaces of lower rank.  The group $G(\RR)$ acts on
$\Dstar_\sigma$ and those real boundary components whose normalizers are
defined over $\QQ$ are called the \emph{rational boundary components}%
\footnote{The actual definition is more complicated but is equivalent to
  what is given here under the condition of geometric rationality.}%
.  Under a condition on $\sigma$ now known as \emph{geometric rationality}
\cite{refnCasselmanGeometricRationality}, Satake
\cite{refnSatakeQuotientCompact} constructs a corresponding
compactification $\Xstar_\sigma$ of $X$ by taking the quotient under $\G$
of the union of the rational boundary components (with a suitable
topology).  The compactification $\Xstar_\sigma$ is stratified by
arithmetic quotients of the rational boundary components.

An important example is when $D$ is a Hermitian symmetric space.  In this
case $D$ may be realized as a bounded symmetric domain in $\CC^N$ for some
$N$ and one of the Satake compactifications is homeomorphic to the natural
compactification $\cl{D} \subseteq \CC^N$.  The various real boundary
components are again Hermitian symmetric spaces.  Geometric rationality in
this case was proved by Baily and Borel \cite{refnBailyBorel}; the
resulting compactification $\Xstar$ of $X$ is called the
\emph{Baily-Borel-Satake compactification}.  Baily and Borel prove that
$\Xstar$ has the structure of a normal projective algebraic variety.

More generally consider the case where $D$ is equal-rank.  A \emph{real
equal-rank Satake compactification} is a Satake compactification for which
all real boundary components are equal-rank symmetric spaces.  The possible
real equal-rank Satake compactifications are enumerated in
\cite[(A.2)]{refnZuckerLtwoIHTwo}.  For some equal-rank symmetric spaces,
such as the one associated to $G(\RR)=\SO(4,4)$, such a compactification
does not exist.  On the other hand, if $G(\RR)= \SO(2p,2q+1)$ then a real
equal-rank Satake compactification does exist (and is unique if $p>1$).

In \cite{refnSaperGeometricRationality} we prove that every real equal-rank
Satake compactification is geometrically rational aside from some
$\QQ$-rank $1$ and $2$ exceptions in which $\QQrank G \neq \RRrank G$..
The resulting compactification $\Xstar_\sigma$ is also called a \emph{real
equal-rank Satake compactification}; the Baily-Borel-Satake
compactification is an example.

\begin{thm}[Zucker/Borel Conjecture \protect{\cite{refnLooijenga},
      \cite{refnSaperSternTwo}}]
\label{thmZuckerBorelConjecture}
 Let $\Xstar_\sigma$ be a real equal-rank Satake compactification of an
equal-rank locally symmetric space $X = \G\back D$.  Then there is a
natural quasi-isomorphism $\Ltwo(\Xstar_\sigma;\EE) \cong
\IpC(\Xstar_\sigma;\EE)$ where $p$ is any middle perversity.
\end{thm}

\begin{proof}
(Compare the proof of the Rapoport conjecture in
\cite[\S27]{refnSaperLModules} and the exposition in \cite[\S\S19,
20]{refnSaperCDM}.)  Let $\Xhat$ be the reductive Borel-Serre
compactification and note that $\Ltwo(\Xstar_\sigma;\EE) \cong
R\pi_*\Ltwo(\Xhat;\EE)$, where $\pi\colon \Xhat \to \Xstar_\sigma$ is
Zucker's quotient map \eqref{refnSatakeQuotient}.  For $x$ in a proper
stratum $F$ of $\Xstar_\sigma$, let $\i_x\colon \{x\} \hookrightarrow
\Xstar_\sigma$ denote the inclusion.  By the local characterization of
middle perversity intersection cohomology on a space with even dimensional
strata \cite{refnGoreskyMacPhersonIHTwo}, \cite[V,
4.2]{refnBorelIntersectionCohomology} it suffices to verify
\begin{alignat}{2}
\label{eqnVanishing}
H^i(\i_x^*R\pi_*\Ltwo(\Xhat;\EE)) &= 0,
&\qquad &i \ge (1/2) \codim F, \\
\label{eqnCoVanishing}
H^i(\i_x^!R\pi_*\Ltwo(\Xhat;\EE)) &= 0,
&\qquad &i \le (1/2) \codim F.
\end{alignat}

We consider \eqref{eqnVanishing}; the proof of \eqref{eqnCoVanishing} is
similar.  Let $k\colon \pi^{-1}(x)\hookrightarrow \Xhat$ and observe that
$H(\i_x^*R\pi_*\Ltwo(\Xhat;\EE)) \cong H(k^*\Ltwo(\Xhat;\EE))\cong
H(k^*\Sheaf_{\Xhat}(\Ltwo(E)))$, where for the last step we use Theorem
~\ref{thmLtwoLmodule}.  However one can show%
\footnote{To construct $k^*\Ltwo(E)$ one first restricts $\Ltwo(E)$ to the
locally closed constructible set $\pi^{-1}(F)$ as in
\cite[\S3.4]{refnSaperLModules} and then to the fiber $\pi^{-1}(x)$ as in
\cite[\S3.5]{refnSaperLModules}.}
that $k^*\Sheaf_{\Xhat}(\Ltwo(E)) \cong \Sheaf_{\pi^{-1}(x)}(k^*\Ltwo(E))$
for an $\L$-module $k^*\Ltwo(E)$ on $\pi^{-1}(x)$ (which is itself the
reductive Borel-Serre compactification of a locally symmetric space
\cite[(3.8), (3.10)]{refnZuckerSatakeCompactifications}).  Thus we are
reduced to the corresponding vanishing of $H(k^*\Ltwo(E))$.  However since
$\mS(\Ltwo(E))=\{E\}$ by Corollary ~\ref{corLtwoMicropurity} and since
all real boundary components are equal-rank, we can use \cite[Corollary
~26.2]{refnSaperLModules} to estimate that $d(k^*\Ltwo(E)) < (1/2) \codim
F$ where $d(\M)$ is defined in \eqref{eqnUpperBound}.  Now apply Theorem
~\ref{thmVanishing}.
\end{proof}

\begin{rem*}
Under the weaker hypothesis that all \emph{rational} boundary components
are equal-rank we can still obtain an estimate of $\mS(k^*\Ltwo(E))$
\cite[Proposition ~23.3]{refnSaperLModules}.  However the precise estimate
on $d(k^*\Ltwo(E))$ in the above proof requires the stronger hypothesis
that all \emph{real} boundary components are equal-rank \cite[Corollary
~25.4, Theorem ~26.1]{refnSaperLModules}.
\end{rem*}

\begin{cor}
Under the hypotheses of the theorem, $H_{(2)}(X;\EE) \cong
I_pH(\Xstar_\sigma;\EE)$.
\end{cor}

\begin{proof}
Since $\Ltwo(\Xstar_\sigma;\EE)$) is fine \cite{refnZuckerWarped},
\cite{refnZuckerLtwoIHTwo}, $H(\Ltwo(\Xstar_\sigma;\EE)) \cong
H_{(2)}(X;\EE)$.
\end{proof}


\bibliographystyle{amsplain} \bibliography{borel04}
\end{document}